\documentclass[a4paper,12pt]{article}

\usepackage{eucal, amssymb,amsmath,graphicx, amsfonts, latexsym, verbatim, psfrag}

\renewcommand{\epsilon}{\varepsilon}

\newcommand{\bx}{\mbox{\boldmath$x$}}

\newcommand{\bX}{\mbox{\boldmath$X$}}

\newcommand{\bpi}{\boldsymbol{\pi}}

\newcommand{\thetahat}{{\hat{\theta}}}

\newcommand{\xtilde}{{\tilde{x}}}

\newcommand{\etal}{{\it et al.}}

\newcommand{\nmax}{n_{\mathrm{max}}}

\newcommand{\dateenglish}{\renewcommand*{\today}{%
\number\day \ifcase\day \or
st\or nd\or rd\or th\or th\or th\or th\or th\or th\or th\or
th\or th\or th\or th\or th\or th\or th\or th\or th\or th\or
st\or nd\or rd\or th\or th\or th\or th\or th\or th\or th\or
st\fi\space \ifcase\month \or
January\or February\or March\or April\or May\or June\or July\or
August\or September\or October\or November\or December\fi \space\number\year}}

\numberwithin{equation}{section}

\newcommand{\hfigwidth}{6.75cm}

\begin{document}
\title{Household epidemic models with varying infection response}

%START JMB title material
%\author{Frank Ball \and Tom Britton \and David Sirl} %\authorrunning{Ball \etal}
%\institute{F.~G.~Ball \at School of Mathematical Sciences, University of Nottingham, University Park, Nottingham NG7 2RD, UK. \email{frank.ball@nottingham.ac.uk} \and
%T.~Britton \at Department of Mathematics, Stockholm University, SE-106 91 Stockholm, Sweden. \email{tomb@math.su.se} \and
%D.~J.~Sirl \at School of Mathematical Sciences, University of Nottingham, University Park, Nottingham NG7 2RD, UK. \email{david.sirl@nottingham.ac.uk}}
%END JMB title material

%START 'normal' (preprint) title material
\title{Household epidemic models with varying infection response \footnote{Submitted to Journal of Mathematical Biology.}}
\author{Frank Ball\thanks{School of Mathematical Sciences, University of Nottingham, University Park, Nottingham NG7 2RD, UK.},
Tom Britton\thanks{Department of Mathematics, Stockholm University, SE-106 91 Stockholm, Sweden.},
David Sirl$^\dagger$}
\date{10th May 2010}
%END 'normal' title material

\maketitle

\begin{abstract}
This paper is concerned with SIR (susceptible $\to$ infected $\to$ removed) household epidemic models in which the infection response may be either mild or severe, with the type of response also affecting the infectiousness of an individual. Two different models are analysed. In the first model, the infection status of an individual is predetermined, perhaps due to partial immunity, and in the second, the infection status of an individual depends on the infection status of its infector and on whether the individual was infected by a within- or between-household contact. The first scenario may be modelled using a multitype household epidemic model, and the second scenario by a model we denote by the infector-dependent-severity household epidemic model. Large population results of the two models are derived, with the focus being on the distribution of the total numbers of mild and severe cases in a typical household, of any given size, in the event that the epidemic becomes established.  The aim of the paper is to investigate whether it is possible to determine which of the two underlying explanations is causing the varying response when given final size household outbreak data containing mild and severe cases. We conduct numerical studies which show that, given data on sufficiently many households, it is generally possible to discriminate between the two models by comparing the Kullback-Leibler divergence for the two fitted models to these data.

%\keywords{Household epidemic model \and infector dependent severity \and Kullback-Leibler divergence \and multitype epidemic \and varying infection response \and final outcome data}
%\subclass{92D30 \and 62M05}
\end{abstract}


%\paragraph{MSC classifications:} 

\section{Introduction}

The present paper concerns models for infectious diseases in a community of households, in which the response to disease varies between individuals; we restrict our attention to having two different responses, denoted mild and severe. One common such situation is for example where there are asymptomatic cases showing no or hardly any symptoms but still contributing to the further spread of the disease.

The reason why individuals show different symptoms may vary for different diseases. In the present paper we focus on two potential explanations. The first explanation is that the disease response is determined by individual characteristics, for example someone having partial immunity might become asymptomatic if infected (e.g.~Staalsoe and Hviid (1998) for malaria and Leroy \etal\ (2001) in the context of ebola). The second explanation we consider is where the response depends on the type of infectious contact and/or whom the individual was contacted by. Examples where this seems to be the case are dengue fever (Mangada and Igarashi 1998), measles (Morley and Aaby 1997) and varicella (Mehta and Chatterjee 2010). Ball and Becker (2006) consider the evaluation of vaccination strategies for a model in which infectious cases may be either mild or severe. However, their analysis is based on post-vaccination reproduction numbers rather than on mechanistic models such as those considered in this paper.

The first explanation, where the response is determined by individual characteristics of the infected person, is suitably modelled using a multitype epidemic household (MT-HH) model (Ball and Lyne 2001). In the MT-HH model individuals are categorized into different types; the type of an individual may affect susceptibility to the disease and also response, in particular infectivity, in the event when the individual becomes infected. Quite often an individual's characteristics would not be known, which implies that the proportions of individuals of the different types in the community are unobserved.

The second explanation can be modelled by extending the so-called infector-dependent-severity (IDS) epidemic model of Ball and Britton (2007) to an infector-dependent-severity household  (IDS-HH) epidemic model. In the IDS model, the probability of an individual becoming a mild/severe case depends on the disease response of the person who caused that individual's infection. In the IDS-HH model this may also depend on whether the contact causing the infection was a within- or between-household contact; for example, within-household transmission might have a higher risk of leading to severe infection.

Once the final size distribution of the IDS-HH epidemic model is obtained, together with known such results for the MT-HH model, it is possible to compare the two distributions. The motivation for the present paper lies in this comparison. In particular we pose (and answer) the following question: can final size data from an epidemic outbreak with varying disease response be used to discriminate between the two candidate explanations for why infection response varies?  Except in a few degenerate cases, the answer to the question is yes. This is illustrated numerically by showing that, in the limit as the population size tends to infinity in an appropriate manner, a possible outcome for either of the model is (usually) inconsistent with the other model. This is done by generating ``data'' from the MT-HH (IDS-HH) model and showing that the Kullback-Leibler divergence of the estimated outbreak probabilities from the ``data'' is much smaller when inference is based on the MT-HH (IDS-HH) model than when it is based on the IDS-HH (MT-HH) model. The final size outcome probabilities for the IDS-HH epidemic are obtained by numerically solving a set of differential equations, and the final size outcome probabilities for the MT-HH model are obtained numerically by solving a set of balance equations. Consequently, we have no analytical results ``proving'' that the two models are inconsistent -- our arguments are instead based on numerical studies. We also consider data generated from finite populations and use a simulation study to demonstrate that it is possible to discriminate between the models using a pseudolikelihood approach.

The paper is organised as follows. In Section~\ref{sec-MT-HH} we define the MT-HH model and review final outcome results for that model. In Section~\ref{sec-IDS-HH} we define the IDS-HH  model and derive an appropriate determinstic approximation to it.  In Section~\ref{sec:FSproperties} we compare and contrast the final size outcomes of the two models via simulation studies, which (i) confirm the applicability of the asymptotic results to finite populations, (ii) strongly suggest that, as proved by Ball and Lyne 2010 for the MT-HH model, the final outcome of the IDS-HH model satisfies a central limit theorem and (iii) shed light on some interesting differences between the models. In Section~\ref{sec:modelDisc} we show numerically that inference from final outbreak data makes it possible to distinguish between the two models, using both infinite populations, as described in the previous paragraph, and also finite populations, where a pseudolikelihood approach (cf.~Ball and Lyne 2010) is applied to simulated data.  The paper ends with a short discussion in Section~\ref{disc}.

\section{The multitype household model}
\label{sec-MT-HH}

\subsection{Model definition}
\label{sec:MTdef}

The multitype household epidemic was first analysed in depth by Ball and Lyne (2001), see also Becker and Hall (1996) and Britton and Becker (2000). We now describe this model using slightly different notation.  With the present application in mind, we restrict the model to two types and exponentially distributed infectious periods.

Each individual is categorized as being a mild or a severe type, with the interpretation that if infected, the individual will become this type of infective. Additional to this, individuals reside in households. Let $m_{k,s}$ denote the number of households having $k$ mild and $s$ severe individuals, let $m_n=\sum_{k=0}^nm_{k,n-k}$ denote the number of households of size $n$, and let $m=\sum_{n=1}^\infty m_n (= \sum_{k,s} m_{k,s})$ denote the total number of households.  Further, let $N=\sum_{n=1}^\infty n m_n$ denote the total population size, which is assumed to be finite. Our analysis is of the limiting situation in which the total number of households $m$ tends to infinity in such a way that $m_n/m \to \rho_n$ ($n=1,2,\ldots$), where $\sum_{n=1}^\infty \rho_n =1$ and the limiting mean household size $\mu_H=\sum_{n=1}^\infty n\rho_n$ is finite. It would rarely be the case that the type of an individual is known, so we assume that each individual is a mild case with probability $\beta_M$ (and severe with probability $1-\beta_M$), with the types of different individuals being mutually independent. It then follows that the number of mild cases in a household of size $n$ is binomially distributed. Hence, in a large community, we have that $m_{k,s}/m_{k+s}\approx \binom{k+s}{k}\beta_M^k(1-\beta_M)^{s}$, and this holds with equality in the limiting situation described above.

The disease spreads according to the following rules. Initially, a small given number of individuals are infected (from some external force) and the remaining individuals are susceptible. During his/her infectious period a mild infectious individual has (global) infectious contacts with any given other mild individual at rate $\lambda_{MM}^{(G)}/N$ and with any given severe individual at rate $\lambda_{MS}^{(G)}/N$. Similarly, an infectious severe individual has (global) infectious contacts with any given mild individual at rate $\lambda_{SM}^{(G)}/N$ and with any given other severe individual at rate $\lambda_{SS}^{(G)}/N$. Additionally, an infectious mild individual has (local) infectious contact with any given other mild member of his/her household at rate $\lambda_{MM}^{(L)}$ and with any given severe mamber of his/her household at rate $\lambda_{MS}^{(L)}$. The corresponding rates for local infectious contacts of an infectious severe individual are $\lambda_{SM}^{(L)}$ and  $\lambda_{SS}^{(L)}$. An `infectious contact' is defined as a contact which results in infection if the other individual is susceptible -- otherwise the contact has no effect. The infectious period of all individuals follow  exponential distributions, with rates $\gamma_M$ and $\gamma_S$ for mild and severe infectives, respectively.  All contact processes are governed by homogeneous Poisson processes, having rates as above.  Further, all infectious periods and all contacts processes (whether or not either or both of the indivdiduals involved are the same) are assumed to be mutually independent.  We assume that infected individuals are able to infect other individuals as soon as they have become infected, i.e.~there is no latent period.  Once an inidividual's infectious period is over, he/she recovers and becomes immune to further infection.  The absence of a latent period, though unrealistic for most, if not all, human diseases, has no consequence for our present purpose, since the distribution of the final outcome of the MT-HH model is not changed if an almost surely finite latent period is incorporated (provided the rest of the model is the same). 

The MT-HH model has the following 11 parameters:
%JMB version
%$\theta^{(MT)}=(\lambda_{MM}^{(G)}, \lambda_{MS}^{(G)},\newline \lambda_{SM}^{(G)}, \lambda_{SS}^{(G)},  \lambda_{MM}^{(L)}, \lambda_{MS}^{(L)}, \lambda_{SM}^{(L)}, \lambda_{SS}^{(L)}, \gamma_M, \gamma_S, \beta_M)$.
%Normal version
$\theta^{(MT)}=(\lambda_{MM}^{(G)}, \lambda_{MS}^{(G)}, \lambda_{SM}^{(G)}, \lambda_{SS}^{(G)}, \lambda_{MM}^{(L)}, \newline \lambda_{MS}^{(L)}, \lambda_{SM}^{(L)}, \lambda_{SS}^{(L)}, \gamma_M, \gamma_S, \beta_M)$.
 Later we consider final size data for this model. In that situation we can, and hence do, assume without loss of generality that $\gamma_M=\gamma_S=1$. (The final outcome of a closed-population stochastic SIR epidemic of this type can be obtained by considering a random directed graph whose vertices are the individuals in the population and, for any vertices  $i \ne j$, there is a directed edge from $i$ to $j$ if and only if individual $i$, if infected, has infectious contact with individual $j$; see, for example, Pellis \etal\ (2008). The set of people who are ultimately infected by the epidemic is given by the those individuals for which there is a chain of directed edges leading to them {\em from} an initial infective. Thus if, for example, $\gamma_M\ne 1$, we can divide all infection rates {\em from} mild infectives by $\gamma_M$ and then set $\gamma_M=1$ without changing the probability measure of the above random directed graph, and hence without changing the final outcome distribution.  Note that this directed random graph explains also the above comment concerning a latent period.) It is also shown in Ball \etal\ (2004) that the 4 global infection parameters are not uniquely identifiable from final size data -- what is identifiable are two separate linear combinations of these four parameters (details are given at the end of Section~\ref{sec:MTLargePop}). In conclusion we hence have 7 parameters that are identifiable from final size data for the MT-HH model.

\subsection{Large population properties of the MT-HH model}
\label{sec:MTLargePop}
The MT-HH model is closely related to the model analysed in Ball and Lyne (2001). (The latter model allows for arbitrarily many types, non-random allocation of types of individuals to households and arbitrary but specified infectious period distributions.) Using essentially the same argument as in Ball and Lyne (2001), the MT-HH model possesses a threshold parameter $R_\ast$, a reproduction number for the proliferation of infected households, which determines whether or not an epidemic started with few initial infectives can become established in a large population.  We now consider the final outcome of such an epidemic that becomes established, so implicitly we assume that $R_\ast$ is above its threshold value of 1.  For $n=1,2,\ldots$ and $r_M,r_S=0,1,\ldots$ such that $r_M+r_S\leq n$, let $p_n^{(MT)} (r_M,r_S | \theta^{(MT)})$ denote the limiting fraction of households of size $n$ that have $r_M$ mild cases and $r_S$ severe cases at the end of an epidemic that becomes established, where the limit is as the total number of households $m \to \infty$.  An outline derivation of a method for determining $p_n^{(MT)} (r_M,r_S|\theta^{(MT)})$ is given below.  It is a slight adaptation of the argument used in Ball and Lyne (2001), which should be consulted for further details.

It is fruitful to consider first the following two-type single-household epidemic model proposed by Addy \etal\ (1991).  Suppose that the household is of size $n$, that it contains $k$ mild individuals and $n-k$ severe individuals, and that all $n$ individuals are initially susceptible.  During the course of the epidemic, individuals avoid infection from outside of the household independently, with probabilities $\pi_M$ and $\pi_S$ for mild and severe individuals, respectively.  The local spread within the household is governed by the same disease dynamics as in the MT-HH model. Write
\[
\Lambda^{(L)} = \left[ \begin{array}{ll}
\lambda_{MM}^{(L)}&\lambda_{MS}^{(L)}\\
\lambda_{SM}^{(L)}&\lambda_{SS}^{L}\end{array}\right], \quad
\bpi = (\pi_M,\pi_S)%, \quad
%\bgamma =(\gamma_M,\gamma_S)
\]
and denote this single-household epidemic model by $E^{(n,k)} ( \Lambda^{(L)}, \bpi )$.  (Recall that we assume that $\gamma_M=\gamma_S=1$.)  Let $Z_M^{(n,k)}$ and $Z_S^{(n,k)}$ denote respectively the numbers of mild and severe removed cases in the household at the end of the single-household epidemic, let
$p^{(n,k)} (i,j | \Lambda^{(L)}, \bpi) = P(Z_M^{(n,k)} = i, Z_S^{(n,k)}=j)$ ($0 \leq i \leq k$, $0 \leq j \leq n-k$),
$\mu_M^{(n,k)} ( \Lambda^{(L)}, \bpi ) = E[Z_M^{(n,k)}]$ and 
$\mu_S^{(n,k)} ( \Lambda^{(L)} , \bpi  ) = E[Z_S^{(n,k)}]$.  The probabilities $p^{(n,k)} (i,j | \Lambda^{(L)} , \bpi )$ may be determined using the following triangular system of linear equations (see Addy \etal\ (1991, Equation~(4))):
\begin{align*}
\sum_{i=0}^{i_1} \sum_{j=0}^{j_1} \frac{
\binom{k-i}{i_1-i}
\binom{n-k-j}{j_1-j}
p^{(n,k)} (i,j | \Lambda^{(L)}, \bpi)}
{
\pi_M^{k-i_1} \pi_S^{n-k-j_1} (h_M (i_1,j_1))^i (h_S (i_1,j_1))^j}
& =
\binom{k}{i_1}
\binom{n-k}{j_1} \\
 & (0 \leq i_1 \leq k, \, 0 \leq j_1 \leq n-k),
\end{align*}
where
\[
h_M (i_1,j_1)=\frac{1}{1 + (k-i_1) \lambda_{MM}^{(L)} + (n-k-j_1) \lambda_{MS}^{(L)}}
\]
and
\[
h_S (i_1,j_1)=\frac{1}{1 + (k-i_1) \lambda_{SM}^{(L)} + (n-k-j_1) \lambda_{SS}^{(L)}}.
\]
The means $\mu_M^{(n,k)} (\Lambda^{(L)}, \bpi )$ and
$\mu_S^{(n,k)} ( \Lambda^{(L)} , \bpi  )$ are easily computed once the probabilities
%normal version
\\
$p^{(n,k)} (i,j | \Lambda^{(L)} , \bpi )$ have been obtained.

Returning to the MT-HH model, suppose that there are few initial infectives, and let $z_M$ and $z_S$ denote respectively the proportions of individuals that are ultimately mild removed and severe removed, respectively.  Then, if the total population size $N$ is large, the probability that a given mild susceptible avoids \textit{global} infection throughout the course of the epidemic is approximately
$\pi_M = \exp [-(Nz_M \frac{\lambda_{MM}^{(G)}}{N} + Nz_S \frac{\lambda_{SM}^{(G)}}{N} )]
= \exp [-(z_M \lambda_{MM}^{(G)}  + z_S \lambda_{SM}^{(G)}  )]$.
The corresponding probability for a given severe susceptible is
$\pi_S = \exp [ - (z_M \lambda_{MS}^{(G)}  + z_S \lambda_{SS}^{(G)})]$.
In the limit as the number of households $m \to \infty$, the approximate probabilities $\pi_M$ and $\pi_S$ become exact and whether or not distinct individuals avoid global infection become independent.  It follows that, in the event of an epidemic becoming established, the final outcome within a typical household of size $n$, that initially contained $k$ mild and $n-k$ severe susceptibles is distributed according to the final outcome of the single-household epidemic $E^{(n,k)} ( \Lambda^{(L)}, \bpi )$, with $\bpi=(\pi_M, \pi_S )$ given by
\begin{equation}
\pi_M = \exp [-(z_M \lambda_{MM}^{(G)} + z_S \lambda_{SM}^{(G)})] \text{\quad and\quad}
\pi_S = \exp [-(z_M \lambda_{MS}^{(G)} + z_S \lambda_{SS}^{(G)})].\label{pimpis}
\end{equation}
Note that the expected final number of mild removal cases in a household chosen uniformly at random is given by $\mu_H z_M$.  Thus, by conditioning on first the size of and then the number of mild individuals in such a randomly chosen household, we have that
\begin{equation}
\mu_H z_M  = \sum_{n=1}^{\infty} \rho_n \sum_{k=0}^n \binom{n}{k} \beta_M^k (1-\beta_M)^{n-k} \mu_M^{(n,k)}  (\Lambda^{(L)}, \bpi).\label{muhzm}
\end{equation}
A similar argument shows that
\begin{equation}
\mu_H z_S = \sum_{n=1}^{\infty} \rho_n \sum_{k=0}^n
\binom{n}{k}
\beta_M^k (1-\beta_M)^{n-k} \mu_S^{(n,k)}
(\Lambda^{(L)}, \bpi ).\label{muhzs}
\end{equation}
After substituting for $\bpi$ from (\ref{pimpis}), equations (\ref{muhzm}) and (\ref{muhzs}) give a pair of non-linear equations for
$(z_M,z_S)$.  These equations always have the solution $(z_M , z_S ) = (0,0)$.  If $R_\ast \leq 1$ this is the only solution, whilst if $R_\ast > 1$ there is (subject to very mild conditions on the parameters) a unique second solution in $[0,1]^2$, $(z_M^\ast , z_S^\ast )$ say, giving the proportions of individuals that are ultimately mild and severe removed in the event of an epidemic that becomes established.  It follows that, if $\bpi ^\ast = ( \pi_M^\ast , \pi_S^\ast )$ is obtained by substituting $(z_M, z_S) = (z_M^\ast , z_S^\ast )$ in (\ref{pimpis}), then, for $n=1,2,\ldots$ and $0 \leq r_M+r_S\leq n$,
\begin{align}
p_n^{(MT)} (r_M, r_S | \theta^{(MT)})=\sum_{k=r_M}^{n-r_S} \binom{n}{k} \beta_M^k (1-\beta_M)^{n-k} p^{(n,k)} &(r_M, r_S | \Lambda^{(L)} , \bpi ^\ast ) \label{pnMT}%\\
%&(r_M, r_S \geq 0,\, r_M + r_S \leq n ). \nonumber
\end{align}
Calculating these final size probabilities numerically is relatively straightforward and follows exactly this procedure. Having substituted for $\bpi$ from~\eqref{pimpis}, we first solve (numerically) the balance equations~\eqref{muhzm} and~\eqref{muhzs} to find $(z_M^\ast,z_S^\ast)$, substitute this into~\eqref{pimpis} to find $(\pi_M^\ast,\pi_S^\ast)$, then use~\eqref{pnMT} to calculate the final size distributions $\{p_n^{(MT)} (r_M, r_S | \theta^{(MT)})\}$.

Note from~\eqref{pimpis} that any $(\lambda^{(G)}_{MM}, \lambda^{(G)}_{MS}, \lambda^{(G)}_{SM}, \lambda^{(G)}_{SS})$ satisfying $z_M^\ast \lambda^{(G)}_{MM} + z_S^\ast \lambda^{(G)}_{SM} = -\log \pi_M^\ast$ and $z_M^\ast \lambda^{(G)}_{MS} + z_S^\ast \lambda^{(G)}_{SS} = -\log \pi_S^\ast$ yields the same final size probabilities $\{p_n^{(MT)} (r_M, r_S | \theta^{(MT)})\}$, so only these two linear combinations of $(\lambda^{(G)}_{MM}, \lambda^{(G)}_{MS}, \lambda^{(G)}_{SM}, \lambda^{(G)}_{SS})$ and not the individual global infection rates are identifiable from final size data. Thus when fitting the MT-HH model we estimate $(\pi_M,\pi_S)$ rather than $(\lambda^{(G)}_{MM}, \lambda^{(G)}_{MS}, \lambda^{(G)}_{SM}, \lambda^{(G)}_{SS})$.

\section{The IDS household model}
\label{sec-IDS-HH}

\subsection{Model definition}
\label{sec:IDSdef}

The infector-dependent-severity household model is an epidemic model where infected individuals may, upon infection, become either severely infected or mildly infected, and the probability that an infected individual becomes mildly (or severely) infected may depend on both the type of its infector and whether the infectious contact is local or global. Additionally, individuals reside in households and the transmission rate is typically appreciably higher between individuals sharing a household. The model is defined as follows.

Assume that there are $N$ individuals in total, and that each individual resides in a household. Let $m_n$ denote the number of households of size $n$ and let $m=\sum_{n=1}^\infty m_n$ denote the total number of  households. We consider the limiting situation in which the population size tends to infinity in the same way as described in Section~\ref{sec:MTdef}. Initially there are $k_M^{(m)}$ mild infectives and $k_S^{(m)}$ severe infectives, with the remaining individuals assumed to be susceptible. (The locations of the initial infectives is discussed later.)  Mild infectives recover and become immune at rate $\gamma_M$ and severe infectives recover and become immune at rate $\gamma_S$.  Thus, the infectious periods of infectives are assumed to follow exponential random variables, with parameter depending on whether an infective is a mild or a severe case. While infectious, a mild infective makes global infectious contacts with any given individual at rate $\lambda_M^{(G)}/N$. If a contacted person is susceptible he/she becomes mildly infected with probability $p_{MM}^{(G)}$ and severely infected with probability $1-p_{MM}^{(G)}$; if a contacted person is already infected then the contact has no effect. Additionally, a mild infective has contact with any (other) given household member (local contact) at rate $\lambda_M^{(L)}$, and such a contacted individual, if susceptible, becomes a mild infective with probability $p_{MM}^{(L)}$ and severe infective with probability $1-p_{MM}^{(L)}$. Severe infectives have contacts according to the same rules, although with parameters $\lambda_S^{(G)}/N$, $p_{SM}^{(G)}$, $\lambda_S^{(L)}$ and $p_{SM}^{(L)}$. All contact processes and infectious periods are assumed to be mutually independent. The epidemic continues until there is no (mild or severe) infective present, when the epidemic stops.

The parameters of the IDS-HH model are
%JMB version
%$\theta^{(IDS)}=(\lambda_M^{(G)},\lambda_S^{(G)}, \lambda_M^{(L)}, \lambda_S^{(L)},\newline p_{MM}^{(G)}, p_{SM}^{(G)}, p_{MM}^{(L)}, p_{SM}^{(L)}, \gamma_M, \gamma_S)$.
%Normal version
$\theta^{(IDS)}=(\lambda_M^{(G)},\lambda_S^{(G)}, \lambda_M^{(L)}, \lambda_S^{(L)}, p_{MM}^{(G)}, p_{SM}^{(G)}, p_{MM}^{(L)}, \newline p_{SM}^{(L)}, \gamma_M, \gamma_S)$.
 Note that rescaling time does not change the final outcome of an epidemic, so, without loss of generality, we may assume that e.g.\ $\gamma_M=1$, whence there are 9 parameters that are, in principle, identifiable from final outcome data. Note also that the directed random graph argument used for the final outcome of the MT-HH model fails to hold for the IDS-HH model, since the distribution of edges emanating from any given individual depends on the type of that individual, which is not determined at the outset of the epidemic and indeed depends on the temporal behaviour of the epidemic.  Thus, fixing $\gamma_M$ as well as $\gamma_S$ would involve a loss of generality and the distribution of the final outcome of the IDS-HH model is generally not invariant to a latent period.

\subsection{Large population properties of the IDS-HH model}
\label{sec:IDSLargePop}
Suppose that the epidemic starts at time $t=0$ and for $t \ge 0$, let $X^{(m)}_{n:i,j,k,\ell}(t)$ denote the number of households of size $n$ that at time $t$ have $i$ mild infectives, $j$ severe infectives, $k$ mild removed individuals and $\ell$ severe removed individuals. Assume now that there is a maximal household size $n_{\max}$, so $\rho_n=0$ for all $n>n_{\max}$.  For $t \ge 0$, let $\bX^{(m)}(t)$ be the vector obtained by letting $n,i,j,k,\ell$ vary over all possible feasible values; viz.~$n=1,2,\ldots,n_{\max}, \, 0 \le i+j+k+\ell \le n$.  Then $\{\bX^{(m)}(t):t\ge0\}$ is a density-dependent Markov population process which can be analysed using theory developed in Ethier and Kurtz (1986, Chapter~11).

Suppose that $m^{-1}\bX^{(m)}(0) \to \bx(0)$ as $ m \to \infty$, where $\bx(0)$ satisfies 
%JMB version
%$\sum_{n=1}^{\nmax} \newline \sum_{k,\ell} x_{n:0,0,k,\ell}(0) < 1$,
%Normal version
$\sum_{n=1}^{\nmax} \sum_{k,\ell} x_{n:0,0,k,\ell}(0) \newline < 1$,
 so a strictly positive fraction of the population is initially infected in the limit as $m \to \infty$.  Then the above-mentioned theory of Ethier and Kurtz (1986) shows that the IDS-HH epidemic process scaled by $m$, $\bar \bX (t):=\bX (t)/m$, converges in probability to a vector of deterministic functions defined by a set of differential equations. More precisely, the component $\bar X_{n:i,j,k,\ell}(t)=X_{n:i,j,k,\ell}(t)/m=\rho_n X_{n:i,j,k,\ell}(t)/m_n$ converges to $\rho_n \tilde{x}_{n:i,j,k,\ell}(t)$ defined below. The interpretation of $\tilde {x}_{n:i,j,k,\ell}(t)$ is hence the (asymptotic) fraction of the size-$n$ households that at time $t$ have $i$ mild infectives, $j$ severe infectives, $k$ mild removed individuals and $\ell$ severe removed individuals. Using this notation we can define the (asymptotic) fraction mildly and severely infected
at time $t$, $i_M(t)$ and $i_S(t)$ respectively, by
\begin{align*}
i_M(t)&=\sum_{n,i,j,k,\ell} i\rho_n \tilde{x}_{n:i,j,k,\ell}(t) / \mu_H \\
i_S(t)&=\sum_{n,i,j,k,\ell} j\rho_n \tilde{x}_{n:i,j,k,\ell}(t) / \mu_H .
\end{align*}

The functions $\tilde{x}_{n:i,j,k,\ell}(t)$ are defined by the following set of differential equations:
\begin{align}
\tilde{x}'_{n:i,j,k,\ell}(t) =& \left( \lambda^{(G)}_Mp^{(G)}_{MM} i_M(t) +  \lambda^{(G)}_Sp^{(G)}_{SM}i_S(t) + \lambda^{(L)}_Mp^{(L)}_{MM}(i-1) + \lambda^{(L)}_Sp^{(L)}_{SM}j\right)\nonumber
\\&
\qquad\qquad\times\left(n-(i-1+j+k+\ell)\right) \rho_n \tilde{x}_{n:i-1,j,k,\ell}(t)\nonumber
\\&+
\left( \lambda^{(G)}_Mp^{(G)}_{MS} i_M(t) +  \lambda^{(G)}_Sp^{(G)}_{SS}i_S(t) + \lambda^{(L)}_Mp^{(L)}_{MS}i + \lambda^{(L)}_Sp^{(L)}_{SS}(j-1)\right)\nonumber
\\&
\qquad\qquad\times\left(n-(i+j-1+k+\ell)\right) \rho_n \tilde{x}_{n:i,j-1,k,\ell}(t)\nonumber
\\&+
\gamma_M(i+1)\rho_n \tilde{x}_{n:i+1,j,k-1,\ell}(t)\nonumber
\\&+
\gamma_S(j+1)\rho_n \tilde{x}_{n:i,j+1,k,\ell -1}(t)\nonumber
\\&-
\left( \lambda^{(G)}_M i_M(t) + \lambda^{(G)}_S i_S(t) + \lambda^{(L)}_Mi + \lambda^{(L)}_Sj\right) \left(n-(i+j+k+\ell)\right) \nonumber
\\&
\qquad\qquad\times\rho_n \tilde{x}_{n:i,j,k,\ell}(t) \nonumber
\\&-
(\gamma_Mi+\gamma_Sj)\rho_n \tilde{x}_{n:i+1,j,k-1,\ell}(t),\label{diff-eq}
\end{align}
with initial values given by $x_{n:i,j,k,\ell}(0)=\rho_n \xtilde_{n:i,j,k,\ell}(0)$.

The differential equation~\eqref{diff-eq} applies to all relevant $(n:i,j,k,\ell)$, i.e.\ where each of the indices are non-negative and $i+j+k+\ell\le n$. Vector components `out of bounds', e.g.~where some index is negative, are defined to be 0, for example $\xtilde_{3:0,-1, 1,1}(t) \equiv 0$.
The first four terms in~\eqref{diff-eq} are for households entering the state $(n:i,j,k,\ell)$, explaining why they have a plus sign. The first term is for a household presently in state $(n:i-1,j,k,\ell)$ having a mild infection and gives the overall rate for such an event to occur. The second term is for $(n:i,j-1,k,\ell)$-households having another severe infection, the third term is for $(n:i+1,j,k-1,\ell)$-households having a mild removal and the fourth term is for the severe removals. The remaining terms describe events that cause a household to leave the state $(n:i,j,k,\ell)$. The fifth term is the overall rate at which susceptibles in a $(n:i,j,k,\ell)$-households become infected, either mildly or severely and the last term is the overall rate at which infectives (mild and severe) in  $(n:i,j,k,\ell)$-households are removed.

Our goal is to obtain $p_n^{(IDS)}(r_M, r_S|\theta^{(IDS)})$, the limiting fraction of households of size $n$ that have $r_M$ mild and $r_S$ severe cases at the end of the epidemic.  If the numbers of initial mild and severe infectives, $k_M^{(m)}$ and $k_S^{(m)}$, are held fixed as $m \to \infty$, then, for large $m$, the epidemic can become established only if the household reproduction number $R_*$ is strictly larger than one.  (The reproduction number $R_*$ can be obtained by approximating the process of infected households by a two-type branching process, the type of an infected household being the type of its initial case; we omit the details as $R_*$ is not required for the present paper.)  Ideally, we would like to be able to calculate  $p_n^{(IDS)}(r_M, r_S|\theta^{(IDS)})$ for an epidemic that becomes established under these conditions. However, if $k_M^{(m)}$ and $k_S^{(m)}$ are held fixed, then $\sum_{n=1}^{\nmax}\sum_{k,\ell} x_{n:0,0,k,\ell}(0)=1$ and the theory of Ethier and Kurtz (1986) cannot be applied directly. Thus we assume instead that a very small, but strictly positive, fraction of individuals are initially infected and approximate $p_n^{(IDS)}(r_M, r_S|\theta^{(IDS)})$, by solving the differential equations~\eqref{diff-eq} numerically up to a time when the remaining fraction of infective individuals is negligible. More specifically, we assume that a fraction $f_S=10^{-5}$ of the population is initially severely infective, with these infective individuals being chosen uniformly at random, so
\begin{equation*}
\xtilde_{n:i,j,k,\ell}(0) = \begin{cases}
\binom{n}{j} f_S^j (1-f_S)^{n-j} & \mbox{ if $i=k=\ell=0$}, \\
0 & \mbox{ otherwise.}
\end{cases}
\end{equation*}
We stop the numerical integration at the first time $t'$ when the proportion of the population that is infective, i.e.\ $i_M(t')+i_S(t')$, is less than $\delta=10^{-7}$ ($\ll f_S$). The final size probabilities are then given by $p^{(IDS)}_n(r_M,r_S|\theta^{(IDS)}) = \xtilde_{n:0,0,r_M,r_S}(t')$. However, note that the final size probabilities are essentially insensitive to the initial conditions, provided that the proportion of index cases is sufficiently small.

The theory of Ethier and Kurtz (1986, Chapter~11) can also be used to show that, in the limit as the number of households $m \to \infty$, the fluctuations of the stochastic model $\bX (t)$ about its deterministic limit $\bx (t)$ (defined by $x_{n:i,j,k,\ell}(t)=\rho_n \xtilde_{n:i,j,k,\ell}(t)$), after being suitably scaled, converge to a zero-mean Gaussian process, whose covariance function can, in principle, be determined.  As in Ball and Britton (2007, 2009), this central limit theorem can be extended heuristically to hold also for the end of the epidemic, the time of which tends to infinity as $m \to \infty$, by making a random time scale transformation in which the clock runs at rate $m( \lambda_M^{(G)} I_M (t) + \lambda_S^{(G)} I_S (t))^{-1}$, where $I_M (t)$ and $I_S (t)$ are respectively the total number of mild and severe infectives present at time $t$ in the untransformed process, cf.\ Ethier and Kurtz (1986, pp.~466--467).  This yields a multivariate central limit theorem for the quantities $Z_n (r_M, r_S )$ ($n \geq 1$, $r_M , r_S \geq 0$, $r_M + r_S \leq n$), where $Z_n (r_M, r_S )$ is the number of households of size $n$ which ultimately have $r_M$ mild removed and $r_S$ severe removed individuals.  In principle, it is possible to compute the covariance matrix of the limiting normal distribution numerically, though in practice the required computations are prohibitive except for populations comprising only very small households.  For a population with households of sizes $1,2,\ldots, n_{\max}$, determining the deterministic limit $\bx (t)$ requires solving a system of $n_E^{\max} = \binom{n_{\max}+5}{5}-n_{\max}-1$ differential equations and determining the above covariance matrix requires solving a system of  $\binom{n_E^{\max} + 1}{2}$ differential equations. For $n_{\max}=1,2,3,4,5$, $n_E^{\max} = 4,18,52,121,246$, so while it is perfectly feasible to solve for $\bx(t)$ numerically, that may not be the case for the covariance matrix.

\section{Numerical illustrations of model behaviour}
\label{sec:FSproperties}
To illustrate the asymptotic results given in the previous sections and explore some of the properties of the two models we have presented, we performed simulation studies of both models and compared some of their final size properties. In order to do this we first needed to select values for the parameters of our models.

First we address these parameter choices in the MT-HH model. As mentioned in Section~\ref{sec:MTdef}, the removal rates can without loss of generality be set to unity: $\gamma_M=\gamma_S=1$. The fraction of mild types in the community was set to $\beta_M=0.4$. This value was chosen so that approximately one third of all infected are mild cases (reported by Carrat \etal\ (2008) to be the case for asymptomatic cases regarding influenza). The global contact rates were chosen as $\lambda_{MM}^{(G)}=0.25$, $\lambda_{MS}^{(G)}=0.8$, $\lambda_{SM}^{(G)}=0.8$ and $\lambda_{MM}^{(G)}=1.5$, so severe infectives are more infectious and also mild infectives rarely globally infect mild susceptibles. The corresponding local contact rates were chosen as $\lambda_{MM}^{(L)}=0.2$, $\lambda_{MS}^{(L)}=0.4$, $\lambda_{SM}^{(L)}=0.4$ and $\lambda_{MM}^{(L)}=0.8$, so severe infectives are also more infectious locally and mild infectives have a relatively higher probability of infecting mild susceptibles when compared with global contacts. The absolute values of the two contact matrices were chosen so that approximately 50\% of the population becomes infected, this being a realistic value for influenza (see Ferguson \etal\ (2005)). The relative magnitude of the global and local contact rates was chosen so that both types of contact play a significant role in the spread of infection.

The parameters of the IDS-HH model were chosen to be $\lambda_M^{(G)}=1$, $\lambda_S^{(G)}=2$, $p_{MM}^{(G)}=0.8$, $p_{SM}^{(G)}=0.2$, $\lambda_M^{(L)}=0.5$, $\lambda_S^{(L)}=1$, $p_{MM}^{(G)}=0.5$, $p_{SM}^{(G)}=0.1$, $\gamma_M=1$ and $\gamma_S=2$. These parameter values were chosen for the same reasons as the parameters for the MT-HH model, with the addition that the infectious period was set to be shorter (on average) for severe cases than for mild cases, having in mind asymptomatic individuals who are less likely to `self-quarrantine' as they are unaware of their infection.

The parameter common to both models is the distribution of household sizes. In this paper we consider two different population structures. The first is for the case where households of size 1, 2 and 3 are equally likely and no larger households exist, i.e.\ $\rho_1=\rho_2=\rho_3=1/3$, and is chosen largely for computational convenience. The second population structure corresponds to the household structure of UK in 2003 (found by typing `household sizes' into the search box at http://www.statistics.gov.uk/census2001/census2001.asp), with the simplification that households of size 5 and larger were truncated and all assumed to have size 5 (only 2\% of the households had larger household size than 5, so this truncation should have negligible effect). The household structure for this case is given by $\rho_1=0.29,\rho_2=0.35, \rho_3=0.15, \rho_4=0.14, \rho_5=0.07$. For future reference we denote these distributions by $\rho^{(3)}=(1,1,1)/3$ and $\rho^{(5)}=(29,35,15,14,7)/100$. We also note that 3 is the minimum value of $\nmax$ for both models to be, in principle, identifiable.

For both models we ran 10,000 simulations of systems with 10,000 households, the sizes being given by $\rho=\rho^{(5)}$. (We treated $\rho$ as giving the proportions of households of different sizes, rather than having random household sizes with distribution $\rho$.) In order that minor outbreaks are unlikely, we initiated the epidemics with 10 infectives, each randomly chosen in different households of size 5; in the MT-HH model these individuals may be mild-type or severe-type (with respective probabilities $\beta_M$ and $1-\beta_M$) and in the IDS-HH model we specified that they are all severe cases. For simulations that result in more than 0.15 of the population becoming infected we then recorded the overall final size amongst initial susceptibles and the within-household final sizes amongst households that are initially completely susceptible. (Inspection of histograms (not shown) of the final proportion of individuals infected suggests that this cutoff is appropriate for separating minor and major outbreaks.)

\begin{figure}[!h]
\begin{center}
\includegraphics[width=\hfigwidth]{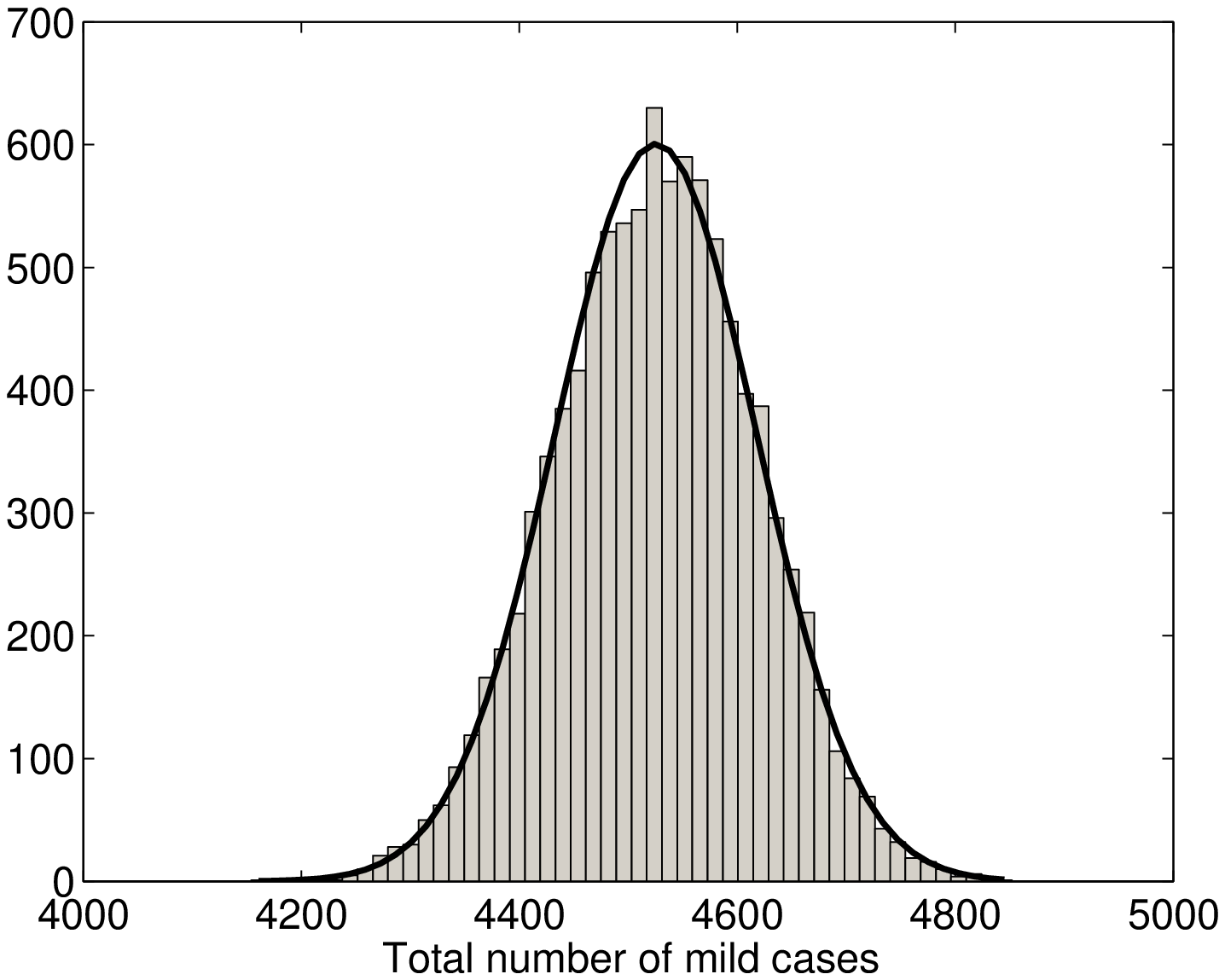}
\includegraphics[width=\hfigwidth]{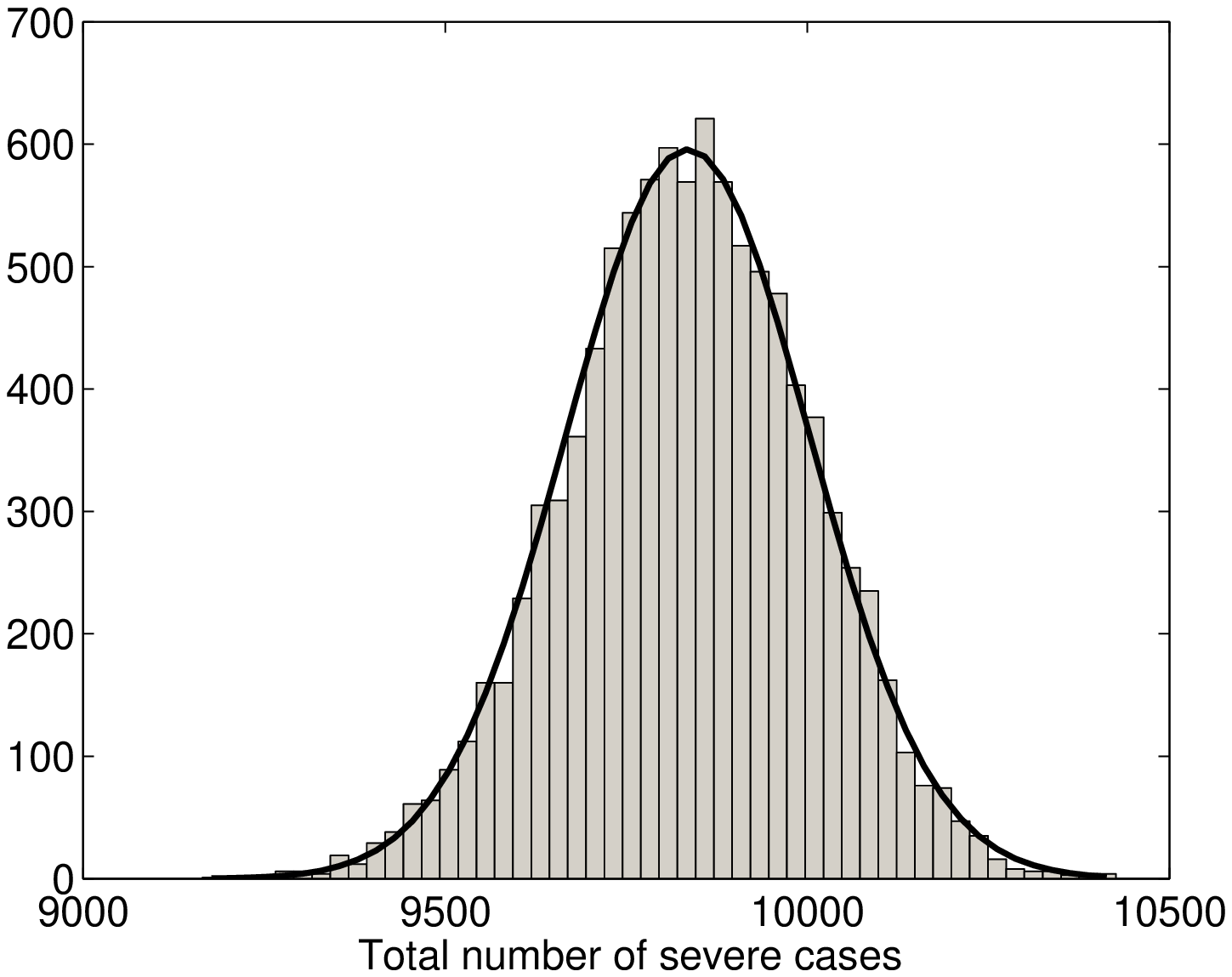}
\caption{Histograms of final outcome of major outbreaks in simulations of the MT-HH model in a community of 10,000 households, with matched normal approximations superimposed.}
\label{fig:MTsimHist}
\end{center}
\end{figure}

\begin{figure}[!h]
\begin{center}
\includegraphics[width=\hfigwidth]{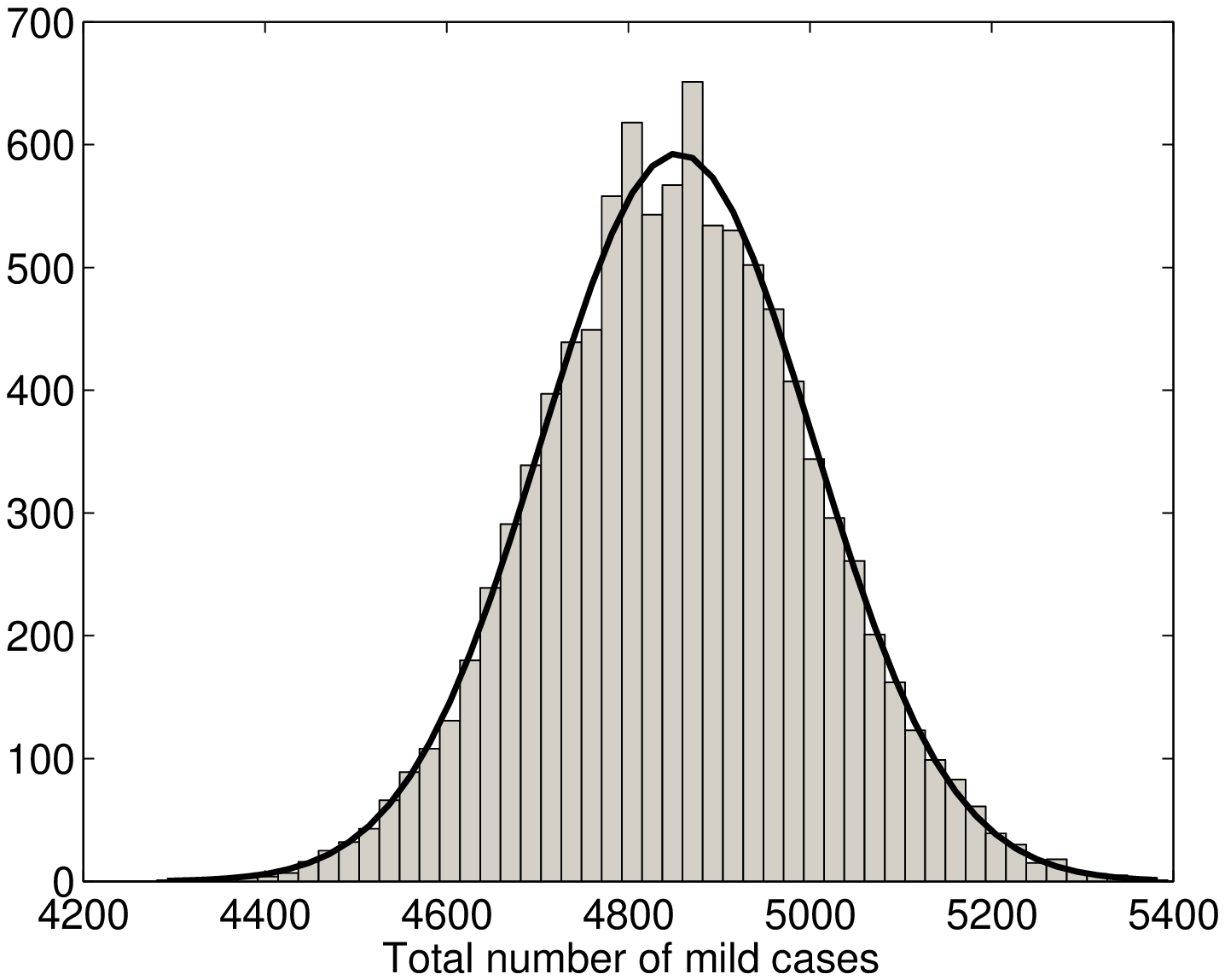}
\includegraphics[width=\hfigwidth]{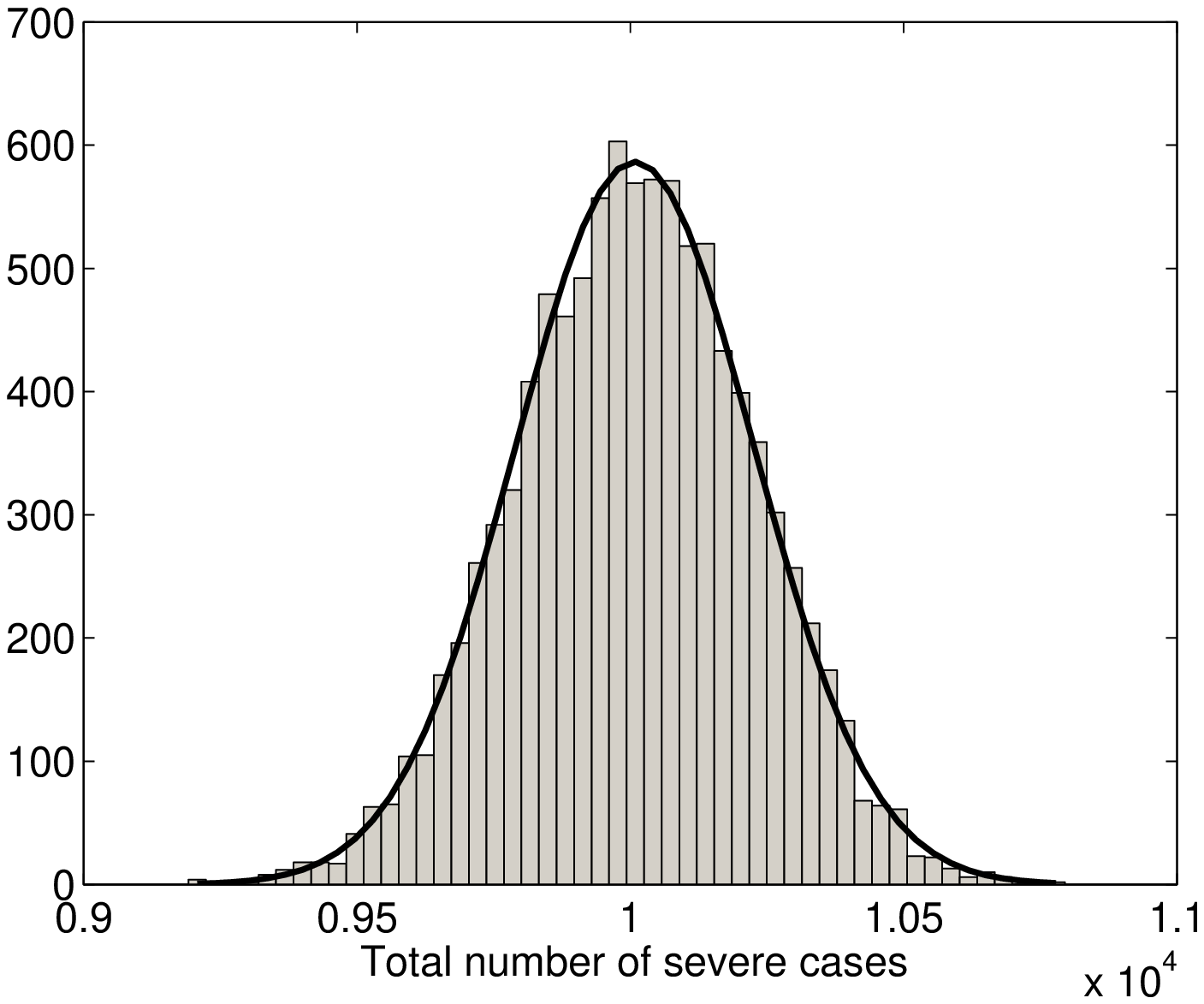}
\caption{Histograms of final outcome of major outbreaks in simulations of the IDS-HH model in a community of 10,000 households, with matched normal approximations superimposed.}
\label{fig:IDSsimHist}
\end{center}
\end{figure}

Figures~\ref{fig:MTsimHist} and~\ref{fig:IDSsimHist} show histograms of the numbers of individuals ultimately mildly and severely infected in, respectively, the 9,992 simulations of the MT-HH model and the 9,993 simulations of the IDS-HH model that resulted in major outbreaks. Overlaid on these histograms are probability density functions (scaled so as the area under them matches that of the histograms) of normal distributions with the same mean and variance. The excellent agreement between the histograms and density functions in Figure~\ref{fig:MTsimHist} is expected in view of the central limit theorem for the MT-HH model of Ball and Lyne (2001) and in Figure~\ref{fig:IDSsimHist} this lends credence to the central limit theorem discussed in Section~\ref{sec:IDSLargePop} above for the final outcome of the IDS-HH model. Though the mean values of these distributions are similar ($\mu_M^{(MT)} \approx 4,\!525$, $\mu_M^{(IDS)} \approx 4,\!854$; $\mu_S^{(MT)} \approx 9,\!835$, $\mu_S^{(IDS)} \approx 10,\!008$) -- indeed the parameter values were chosen with this intention -- it is interesting to note that the variability is rather different in the two models. The spread of the distribution of the number of mild cases in the IDS-HH model is appreciably larger than that in the MT-HH model ($\sigma_M^{(MT)} \approx 93$ and $\sigma_M^{(IDS)} \approx 150$) and, though not to the same extent, the distribution of the number of severe cases is also slightly more spread in the IDS-HH model ($\sigma_S^{(MT)} \approx 167$ and $\sigma_S^{(IDS)} \approx 218$). Part of the reason for this is that in the MT-HH model the types of individuals are determined in advance, but in the IDS-HH model the types of the infected individuals depend on the evolution of the epidemic and some feedback may occur (though with different parameters it might potentially be positive or negative).

Tables~\ref{tab:MTestEpProp} and~\ref{tab:IDSestEpProp} give further information about the within-household outcomes of major outbreaks in the two models. These tables give, for each household size $n$, estimates from the simulations of the probability that a typical individual in a household of size $n$ is ultimately (i) mildly infected ($\hat{p}_M$), (ii) severely infected ($\hat{p}_S$) and (iii) infected ($\hat{p}_{\textrm{INF}}=\hat{p}_M + \hat{p}_S$), and also of the probability that a case in a household of size $n$ is severe ($\hat{p}_S / \hat{p}_{\textrm{INF}})$. The figures in parentheses are the corresponding infinite population asymptotic quantities obtained from $\{p_n^{(MT)}(r_M,r_S | \theta^{(MT)})\}$ and $\{p_n^{(IDS)}(r_M,r_S | \theta^{(IDS)})\}$, respectively.

\begin{table}
\caption{Properties of MT-HH epidemics that become established.}
\begin{equation*}
\renewcommand{\arraystretch}{1.1}
\begin{array}{l|c|c|c|c}
n & \hat{p}_M & \hat{p}_S & \hat{p}_{\textrm{INF}} & \hat{p}_S / \hat{p}_{\textrm{INF}}\\ \hline
1 & 0.1272~(0.1273) & 0.3255~(0.3256) & 0.4527~(0.4529) & 0.7190~(0.7189) \\
2 & 0.1585~(0.1585) & 0.3751~(0.3753) & 0.5336~(0.5337) & 0.7030~(0.7031) \\
3 & 0.1923~(0.1925) & 0.4228~(0.4229) & 0.6151~(0.6154) & 0.6873~(0.6872) \\
4 & 0.2270~(0.2271) & 0.4656~(0.4658) & 0.6926~(0.6929) & 0.6723~(0.6722) \\
5 & 0.2602~(0.2603) & 0.5020~(0.5021) & 0.7621~(0.7624) & 0.6587~(0.6586)
\end{array}
\renewcommand{\arraystretch}{1}
\end{equation*}
\label{tab:MTestEpProp}
\end{table}

\begin{table}
\caption{Properties of IDS-HH epidemics that become established.}
\begin{equation*}
\renewcommand{\arraystretch}{1.1}
\begin{array}{l|c|c|c|c}
n & \hat{p}_M & \hat{p}_S & \hat{p}_{\textrm{INF}} & \hat{p}_S / \hat{p}_{\textrm{INF}}\\ \hline
1 & 0.1815~(0.1822) & 0.2870~(0.2865) & 0.4685~(0.4687) & 0.6126~(0.6113) \\
2 & 0.1969~(0.1976) & 0.3546~(0.3542) & 0.5515~(0.5517) & 0.6430~(0.6419) \\
3 & 0.2095~(0.2104) & 0.4267~(0.4261) & 0.6362~(0.6364) & 0.6707~(0.6695) \\
4 & 0.2190~(0.2196) & 0.4980~(0.4975) & 0.7169~(0.7171) & 0.6946~(0.6937) \\
5 & 0.2229~(0.2250) & 0.5671~(0.5638) & 0.7901~(0.7888) & 0.7178~(0.7147)
\end{array}
\renewcommand{\arraystretch}{1}
\end{equation*}
\label{tab:IDSestEpProp}
\end{table}

In both cases we observe good agreement between the deterministic and estimated stochastic quantities. Also observe that, in both models, the proportion of individuals infected increases with household size $n$. This a consequence of local spread being greater in larger households. In the MT-HH model the proportion of cases that are severe decreases with $n$, whereas this proportion increases with $n$ in the IDS-HH model. In the IDS-HH model this simply reflects the fact that local infections are very likely to result in severe cases. In the MT-HH model, however, an individual's type is determined in advance rather than by the spread of infection and this results in a `saturation effect' of sorts. Note that the proprtion of cases in households of size 1 that are severe (0.7191) is larger than the proportion of individuals that are of severe type (0.6). In larger households more local spread is expected than in smaller households (as $\lambda^{(L)}_{MM}$, $\lambda^{(L)}_{MS}$, $\lambda^{(L)}_{SM}$ and $\lambda^{(L)}_{SS}$ are independent of household size) and the rates are such that severe types are more likely to be infected. Indeed, in a very large household we would expect everyone to be infected, in which case the proportion of cases that are severe must be equal to the proportion of individuals of severe type. For any household size, the proportion of globally contacted individuals that are severe is 0.7191. However, since local spread increases with households size and there is greater scope for local infection amongst mild types than severe types (as relatively fewer are infected globally) the proportion of cases that are severe must decrease with household size.

Clearly the above phenomena depend on the parameter values chosen in the two models. For example, in either model, simply interchanging the labels of the two types results in the opposite effect of household size on $\hat{p}_S / \hat{p}_{\textrm{INF}}$ being observed.

\section{Model discrimination}
\label{sec:modelDisc}
Suppose data from a population with household structure given by $\{\rho_n\}$ are generated from one of the two models, with some given parameters $\theta^{(MT)}$ or $\theta^{(IDS)}$, as appropriate. An important inference, or discrimination, problem in light of two possible models is then whether it is possible to determine which of the models the data come from. It is hard to give an analytical answer to this question since the final size probabilities are not explicit. We address the question with a numerical investigation. For our purposes, the data are the distributions of within-household final sizes $q=\{q_n(r_M,r_S) ,\, 0\leq r_M+r_S \leq n,\, 0\leq n\leq\nmax\}$. We consider the case where $q$ is the asymptotic ($m\to\infty$) final size distribution derived from one or other of the models using the methods described in Sections~\ref{sec:MTLargePop} and~\ref{sec:IDSLargePop}, in order to determine whether the two models actually produce different final size distributions. We also consider the case where $q$ is derived from stochastic simulations of one or other of the two models (i.e.\ with $m$ finite), to determine whether or not any difference between the two models is sufficiently pronounced to be detectable with a dataset that resembles more closely one available in real life.

In the remainder of this section we describe first, in Section~\ref{sec:dataGen}, how we generate the data that we use, both from infinite and finite populations, then discuss, in Section~\ref{sec:modelFitting} how we fit the models to a given final size distribution. In Section~\ref{sec:disc}, we describe our main findings concerning whether the MT-HH and IDS-HH models can be distinguished on the basis of final size data. In Section~\ref{sec:KLmotivation} we motivate our use of the Kullback-Leibler divergence as a tool for model fitting and discrimination and finally, in Section~\ref{sec:ident}, we discuss identifiability issues that arise in fitting the models to final outcome data.

\subsection{Data generation}
\label{sec:dataGen}

Final size data for an infinite population are generated using the methods described in the previous sections. For the MT-HH model we solve equations (2.1)--(2.4) numerically and for the IDS-HH model we solve the differential equations~\eqref{diff-eq} numerically, as described in Sections~\ref{sec:MTLargePop} and~\ref{sec:IDSLargePop}, respectively. To generate final size data for finite populations we simulate an outcome of the relevant stochastic process according to the model description in Section~\ref{sec:MTdef} or~\ref{sec:IDSdef}, as appropriate, (with 10 initial severe infectives, each in separate households of size $\nmax$, to increase the chance of a major outbreak occuring). If a major outbreak does occur (which we take to be more than 0.15 of the population becoming infected) then we calculate the empirical final size distribution considering only the households in that simulation that had no initial infectives. In either case we denote by $q=\{q_n(r_M,r_S) \}$ the `target' household final size distributions that we try to reproduce from the model we choose to fit to the data.

\subsection{Model fitting}
\label{sec:modelFitting}
We now describe the algorithm we use to fit each model to given final size data $q=\{q_n(r_M,r_S)\}$. The goal is to find parameters $\theta$ of the model we are fitting so that the distance between the final size distributions $\{p_n(r_M, r_S|\theta)\}$, corresponding to $\theta^{(MT)}$ or $\theta^{(IDS)}$, and the `target' final size distributions $q=\{q_n(r_M,r_S)\}$ is as small as possible. We measure this distance using the Kullback-Leibler (K-L) divergence
\begin{equation}
f(\theta) = D_{KL} (q || p(\theta))=\sum_{n=1}^{n_{\max}} \rho_n \sum_{r_M,r_S} q_n (r_M,r_S) \log \left( \frac{q_n (r_M,r_S)}{p_n (r_M,r_S | \theta)} \right).
\label{eq:KLdef}
\end{equation}
The use of the K-L divergence is motivated by its well-known relationship with likelihood-based inferential procedures (see, for example, Bishop \etal\ (1975, pp.~344--348)), which is discussed in more detail in Section~\ref{sec:KLmotivation}. We minimise $f(\theta)$ numerically using Matlab's {\tt fmincon} constrained optimisation routine. This requires selecting a starting point $\theta_0$ for the parameters; we choose these starting values independently at random, the rate parameters from an exponential distribution with mean 1 and proportion/probability parameters uniformly from the interval $(0,1)$. In the case of fitting the IDS-HH model we find that the numerical optimisation is more difficult and that it is beneficial to sample several (we use 20) such possible starting points $\theta_0$ and then start the numerical optimisation routine at the best of these points (i.e.\ that with smallest $f(\theta_0)$), so that the numerical routine is more likely to start at a point in parameter space that is at least moderately compatible with the target final size distributions $q$. We describe this process of choosing a starting point for and then running the optimisation routine as a single `run' of our algorithm (i.e.\ model fitting procedure).

Calculating $p(\theta)=\{p_n (r_M,r_S | \theta)\}$ using the methods described in Section~\ref{sec:MTLargePop} or~\ref{sec:IDSLargePop} is straightforward and in principle evaluating $f(\theta)$ is then trivial as long as $p$ has no zero entries, i.e.\ as long as the parameter vector $\theta$ results in the model being super-critical. However, there are numerical problems that can arise when calculating the K-L divergence as in equation~\eqref{eq:KLdef}. These problems arise due to so-called `catastrophic cancellation' which occurs when using the formula~\eqref{eq:KLdef} if $q$ and $p$ differ only slightly. The terms $q_i \log(q_i/p_i)$ are all small (since $p$ and $q$ are close) but are of differing signs (since sometimes $q_i>p_i$ and sometimes vice-versa), thus when the sum $\sum_i q_i \log(q_i/p_i)$ is close to zero there can be catastrophic cancellation and the calculated sum can be wildly inaccurate. We resolve this by using the Taylor series approximation $s\log(s/t) \approx (s-t)^2/2t$ about $s=t$ (cf.\ Bishop \etal\ (1975, Lemma~14.9-1)), which implies
\begin{equation}
f(\theta) = D_{KL} (q || p(\theta)) \approx \sum_{n=1}^{n_{\max}} \rho_n \sum_{r_M,r_S}  \frac{(q_n (r_M,r_S) - p_n (r_M,r_S | \theta))^2}{2p_n (r_M,r_S | \theta)} .
\label{eq:KLapprox}
\end{equation}
This approximation becomes exact as $p\to q$ so using it when the calculated K-L distance is small gives a good approximation and avoids numerical problems. Numerical experiments comparing the calculated values of $f(\theta)$ using~\eqref{eq:KLdef} and~\eqref{eq:KLapprox} show good agreement, improving as $f(\theta)$ becomes smaller (precisely as expected), but when $f(\theta)$ is less than about $10^{-6}$ we begin to see significant disagreement. Therefore, all our calculations of K-L distance initially use~\eqref{eq:KLdef} but if the result is smaller than $10^{-5}$ we recalculate using~\eqref{eq:KLapprox}.

The random starting values of $\theta^{(MT)}$ and $\theta^{(IDS)}$ in our fitting procedure will sometimes be poor (i.e.\ give large values of $f(\theta)$) and result in the optimisation routine staying in a part of parameter space that gives a very poor fit. Thus, when fitting a model to data we run our algorithm many times over to ensure that as much as possible of the parameter spaces are explored. The number of these runs necessary is somewhat variable; this issue is addressed in Section~\ref{sec:ident}.

Initially we focus simply on the smallest of the K-L distances $f(\thetahat)$ of the model from the data that we find for each combination of dataset and model. In Section~\ref{sec:ident} we explore in more detail the variability of the $f(\thetahat)$ from run to run of our algorithm and also examine the behaviour of the corresponding parameter estimates $\thetahat$.

\subsection{Model discrimination}
\label{sec:disc}

\subsubsection{Infinite data}
\label{sec:infPopn}
To determine whether or not each model is capable of producing the final size distributions generated by the other model we fit a given final size distribution to both models, in the expectation that the correct model will fit appreciably better. We find that the correct model can be made to fit as well as we please by tightening the stopping criteria of the numerical optimisation routine but that there is a definite non-zero lower bound for $f(\thetahat)$ when we fit the wrong model. Further details of this are given in Section~\ref{sec:ident}.

We summarise our findings by way of Table~\ref{tab:2x2Inf}, which shows the best fits obtained from 100 runs of our algorithm (measured by $f(\thetahat)$) obtained when fitting the IDS-HH and MT-HH models to the (asymptotic) final size distributions produced from each of the models (with parameter values as in Section~\ref{sec:FSproperties}) with each of the household size distributions $\rho^{(3)}$ and $\rho^{(5)}$. Table~\ref{tab:2x2Inf} demonstrates the significant differences in fit obtained when fitting the two models to each data set (i.e.\ each column of the table). In Table~\ref{tab:2x2Inf} and the following discussion, `data' refers to the model that generated the given final size distributions we fit to and `model' refers to the model we fit to these data.
\begin{table}
\begin{center}
\caption{Best fits of each model to final size distributions obtained from all four combinations of household size distribution and model.}
\renewcommand{\arraystretch}{1.2}
\begin{equation*}
\begin{array}{rrr@{\,\times\,}lr@{\,\times\,}l|r@{\,\times\,}lr@{\,\times\,}l}
 & & \multicolumn{8}{c}{\mbox{data}} \\
 &  & \multicolumn{4}{c}{\rho=\rho^{(3)}} & \multicolumn{4}{c}{\rho=\rho^{(5)}} \\
 & \hfill \vline & \multicolumn{2}{c}{\mbox{MT-HH}} & \multicolumn{2}{c}{\mbox{IDS-HH}} & \multicolumn{2}{c}{\mbox{MT-HH}} & \multicolumn{2}{c}{\mbox{IDS-HH}} \\ \cline{2-10}
 & \mbox{MT-HH\ } \vline &
3.4 & 10^{-11} & 1.5 & 10^{-3} & 2.0 & 10^{-11} & 6.8 & 10^{-3} \\
\raisebox{-6pt}[0pt]{\rotatebox{90}{model}} 
& \mbox{IDS-HH\ } \vline &
4.7 & 10^{-5} & 8.9 & 10^{-9} & 1.1 & 10^{-4} & 3.3 & 10^{-8}
\end{array}
\end{equation*}
\renewcommand{\arraystretch}{1}
\label{tab:2x2Inf}
\end{center}
\end{table}

It can be seen that using the household size distribution $\rho^{(5)}$ which includes households up to size 5 makes no qualitative difference to these conclusions. However, it is interesting to examine the effect of larger households on the above K-L distances. Table~\ref{tab:KLsummands} shows the contribution to the best final K-L distances in Table~\ref{tab:2x2Inf} from households of each size, which amounts to separating out the summands $\rho_n \sum_{r_M,r_S} q_n (r_M,r_S) \log ( q_n (r_M,r_S) / p_n (r_M,r_S | \theta) )$ in equation~\eqref{eq:KLdef}. The breakdown of the best final K-L distances suggests that an appreciably greater contribution to the K-L distances comes from larger households than would be expected simply based upon the proportions of households of different sizes present. This perhaps suggests that data collection effort might be focused somewhat more on larger households; though of course this depends crucially on our assumption that the same transmission parameters apply in households of all sizes.

\begin{table}
\begin{center}
\caption{Breakdown of contribution to final K-L distances by households of different sizes when $\rho=\rho^{(3)}$.}
\renewcommand{\arraystretch}{1.2}
\begin{equation*}
\begin{array}{c|r@{\,\times\,}l r@{\,\times\,}l r@{\,\times\,}l r@{\,\times\,}l}
 & \multicolumn{2}{c}{\mbox{MT-HH model,}} & \multicolumn{2}{c}{\mbox{IDS-HH model,}} & \multicolumn{2}{c}{\mbox{MT-HH model,}} & \multicolumn{2}{c}{\mbox{IDS-HH model,}} \\
n & \multicolumn{2}{c}{\mbox{MT-HH data}} & \multicolumn{2}{c}{\mbox{IDS-HH data}} & \multicolumn{2}{c}{\mbox{IDS-HH data}} & \multicolumn{2}{c}{\mbox{MT-HH data}} \\ \hline

1 & 4.4 & 10^{-12} & 7.6 & 10^{-10} & 2.0 & 10^{-5} & 2.0 & 10^{-7} \\
2 & 1.1 & 10^{-11} & 1.0 & 10^{-9} & 3.2 & 10^{-5} & 1.1 & 10^{-5} \\
3 & 1.8 & 10^{-11} & 7.1 & 10^{-9} & 1.4 & 10^{-3} & 3.6 & 10^{-5} \\ \hline
\rm{total}
  & 3.4 & 10^{-11} & 8.9 & 10^{-9} & 1.5 & 10^{-3} & 4.7 & 10^{-5}
\end{array}
\end{equation*}
\renewcommand{\arraystretch}{1}
\label{tab:KLsummands}
\end{center}
\end{table}

We see (Table~\ref{tab:2x2Inf}) that the final size distribution generated by each model using somewhat realistic parameter values cannot be captured by the other model. To investigate whether this conclusion holds for the models in general, we need to do this comparison for a range of (supercritical) parameter values. We expect that the fits will be poor except possibly for some degenerate cases where the models can produce the same final size distributions.

To test whether one model can reproduce final size data from the other, we select parameters for one model at random, resampling if the corresponding $R_*\leq1$, and calculate the corresponding final size distribution, then fit the other model to this `data'. When selecting the random model parameters to use, each parameter is chosen independently, rate parameters from an exponential distribution with mean 1 and probability parameters uniformly on $[0,1]$. We then repeat this experiment many times so that we explore a range of parameter combinations of the model from which we derive our data. The final K-L distance (best fit) $f(\thetahat)$ that we report for each paremeter combination is the best fit obtained in 5 runs of our algorithm. (When fitting to final size data that is not exactly reproducible by the model we are fitting we find that the variability of $f(\thetahat)$ between runs is very small and thus 5 runs is more than sufficient to be confident that we have found the best-fitting model; see Section~\ref{sec:ident} for further details.)

\begin{figure}
\begin{center}
\psfrag{IDSmodelMTdata}[][]{\tiny (b) fitting IDS-HH model to MT-HH data}
\psfrag{MTmodelIDSdata}[][]{\tiny (a) fitting MT-HH model to IDS-HH data}
\psfrag{log10fthetahat}[][]{\tiny $\log_{10} f(\thetahat)$}
\includegraphics[width=\hfigwidth]{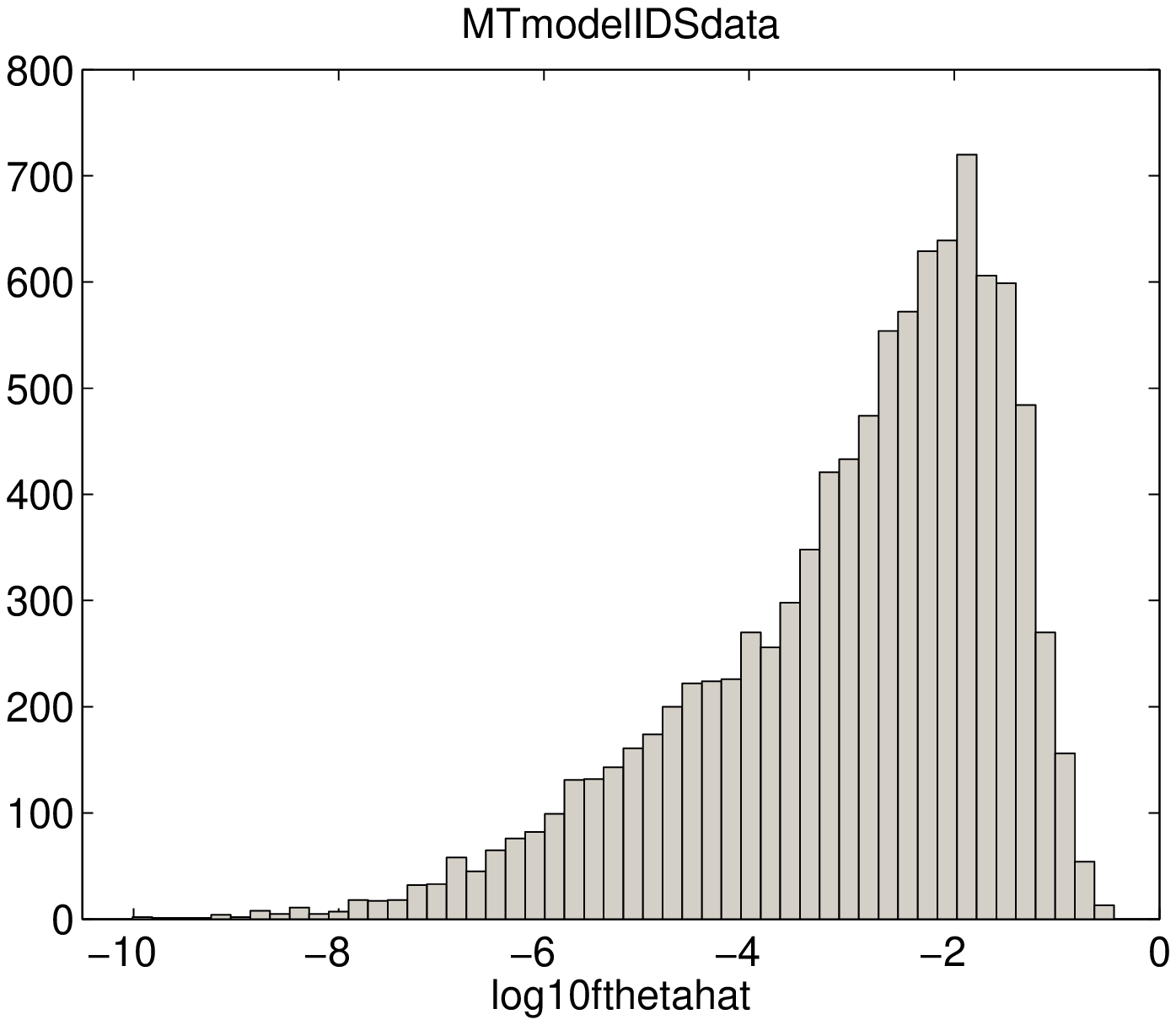}
\includegraphics[width=\hfigwidth]{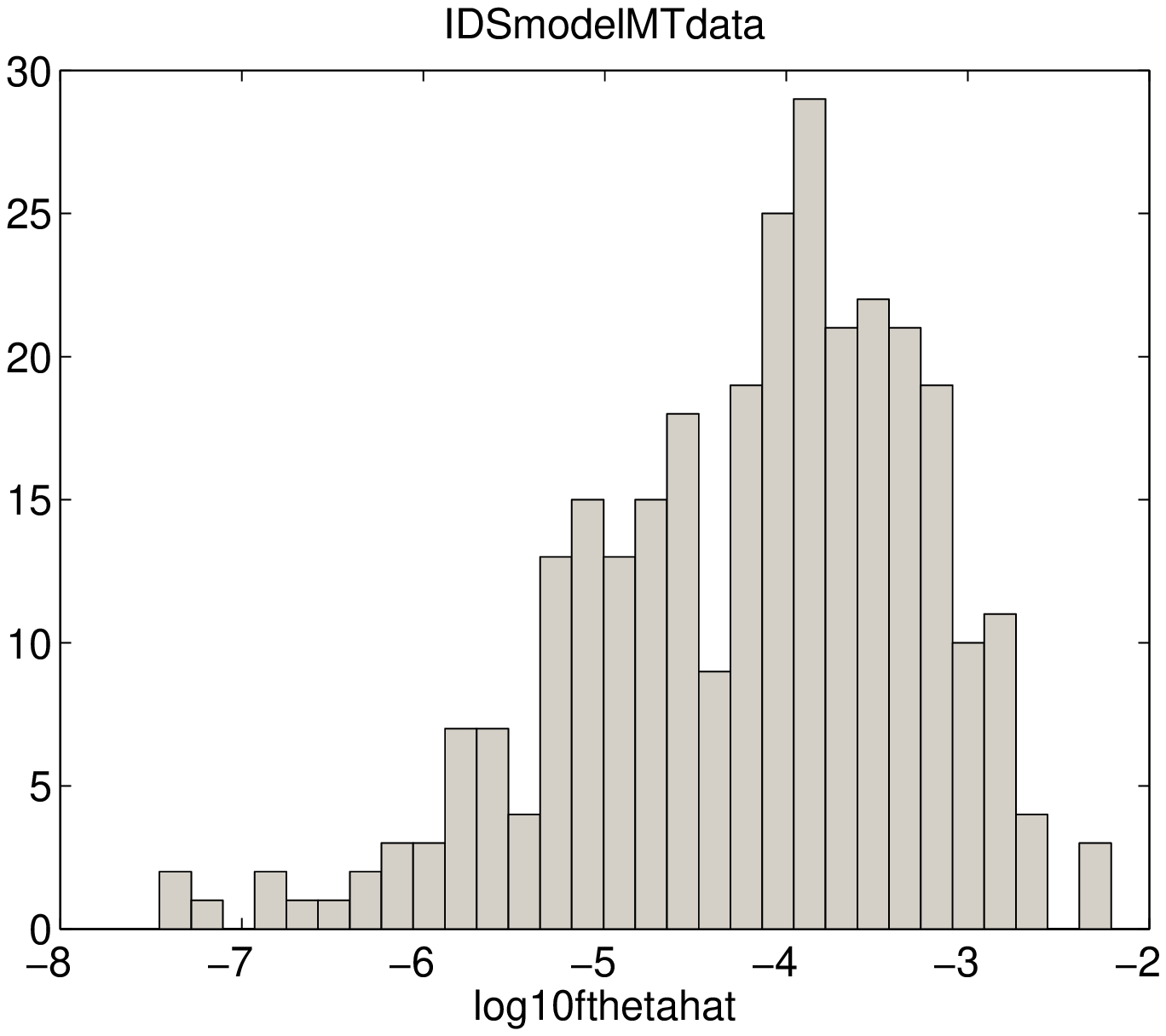}
\caption{Smallest K-L distances $f(\thetahat)$ obtained when fitting the one model to the final size distributions from the other with random (super-critical) parameters.}
\label{fig:InfHist}
\end{center}
\end{figure}

Figure~\ref{fig:InfHist} shows histograms of the best K-L distances $f(\thetahat)$ obtained when fitting one model to final size data generated by the other with random parameters and household size distribution $\rho^{(3)}$. We have fitted the MT-HH model to 10,000 random IDS-HH datasets but owing to the computational expense of fitting the IDS-HH model we have only fitted it to 300 random MT-HH datasets. It can clearly be seen that for most parameter combinations the final size distributions cannot be reproduced by the wrong model. Of course the correct model can reproduce these final size distributions and to confirm this we have also fitted the correct model to many of these data. As expected, the correct model fits appreciably better than the wrong model except in the degenerate cases discussed below when both models can be fitted to the data (details not shown).

%\emph{These calculations/plots are being re-done in light of the problems with small K-L distances. (a) is done and (b) currently has 136 `iterations'. Also need to do some investigation re fitting to correct model so comparison is fair---some preliminary calculations suggest that the true model always fits at least as well.}

Further analysis of the cases where the `wrong' model fits the data relatively well ($f(\thetahat)<10^{-6}$) reveals at least one of the following reasons. In either model, if the process is only just super-critical, i.e.\ $R_*$ is only slightly larger than 1, then many of the quantities $q_n(r_M,r_S)$ in~\eqref{eq:KLdef} are very small so relatively fewer of the summands contribute to the sum and it is somewhat easier for the wrong model to be able to fit the data.

In the IDS-HH model, two further situations arise where the wrong (i.e.\ MT-HH) model fits the final size data quite well: (i) one or both types of individual makes very few local contacts, i.e.\ $\min\{\lambda^{(L)}_M,\lambda^{(L)}_S/\gamma_S\}$ is close to 0 (recall $\gamma_M=1$) and (ii) local contacts by mild and severe individuals are approximately equally likely to cause mild/severe cases, i.e.\ $|p^{(L)}_{MM}-p^{(L)}_{SM}|$ is close to 0. In case (i), within-household spread essentially only involves one type of individual making contacts and the local infection rate and probability can be tuned to produce almost any local final outcome distribution. In case (ii) local infection processes become like those in the MT-HH model because each (locally infected) individual becomes mild or severe (independently) with the same probability $p^{(L)}_{MM} \approx p^{(L)}_{SM}$. In the MT-HH model there are also two further situations where the wrong (i.e.\ IDS-HH) model can reproduce the final outcome data well. The first is if there is essentially only one type of individual, i.e.\ $\beta_M$ is close to 0 or 1; if one type is not present it is trivial that the two models coincide (in the sense that they can produce the same final size distributions). The second case is where the disease is highly globally infectious amongst one type of individual, i.e.\ $\bpi$ has (at least) one element close to~0. Here local transmission is essentially a one-type process and again the models coincide.

\subsubsection{Finite data}
\label{sec:finitePopn}
We have just seen that it is possible to discriminate between the two models using final size data from an infinite population. Real data of course never pertains to an infinite population, so in the present subsection we perform the same type of analysis except that data are now generated from the stochastic models in a community of $m=10,000$ households. The data, generated both from the stochastic MT-HH model and the stochastic IDS-HH model, hence consist of empirical final size distributions rather than the exact asymptotic distributions. The model fitting procedure is exactly the same as before but we now use the empirical final size distribution as the target final size distribution $q$ in~\eqref{eq:KLdef}.

\begin{figure}
\begin{center}
\psfrag{fthetahat}[][]{\tiny $f(\thetahat)$}
\psfrag{titleIDS}[][]{\tiny (a) simulated IDS-HH data}
\psfrag{titleMT}[][]{\tiny (b) simulated MT-HH data}
\includegraphics[width=\hfigwidth]{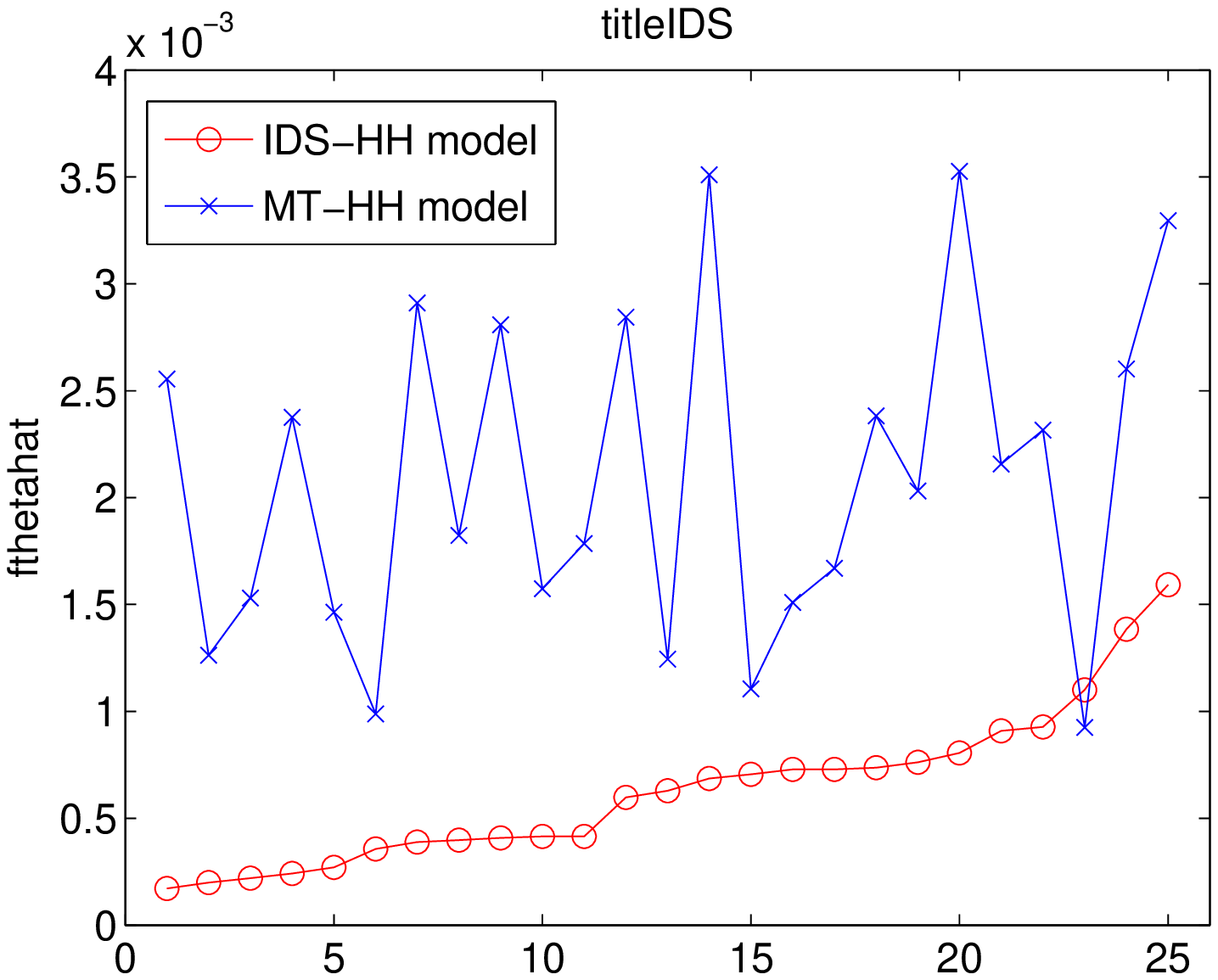}
\includegraphics[width=\hfigwidth]{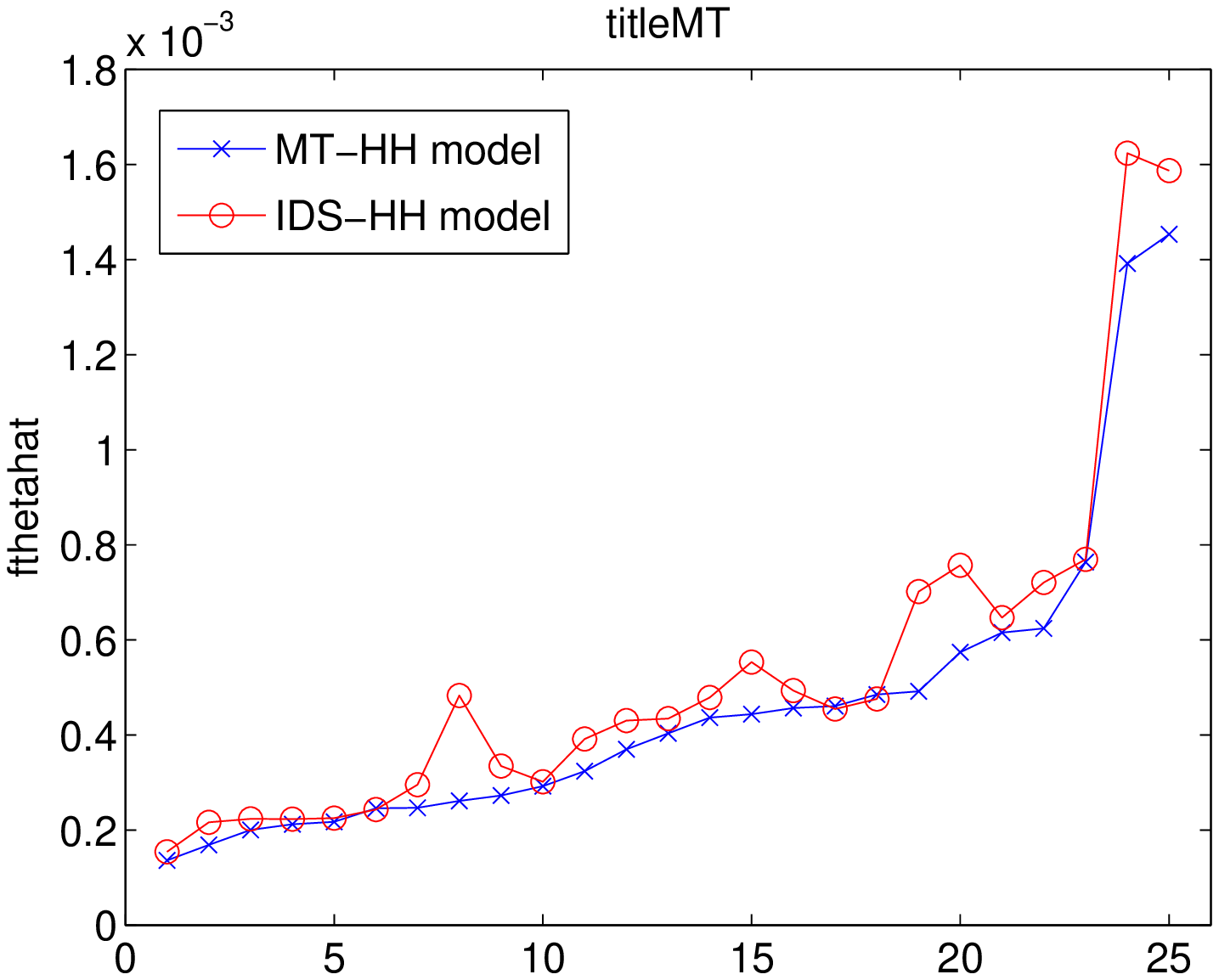}
\caption{Smallest K-L distances $f(\thetahat)$ obtained when fitting both models to the output of 25 separate empirical final size distributions from simulations of each model.}
\label{fig:simGraphs}
\end{center}
\end{figure}

From each model, with parameter valules as in Section~\ref{sec:FSproperties} and household size distribution $\rho^{(3)}$, we generated 25 independent empirical final size distributions taken from simulations on systems of 10,000 households that resulted in a major outbreak and then fitted both models to each empirical final size distribution. Figure~\ref{fig:simGraphs} shows plots of the fit of each model to each dataset, the fit being measured by the smallest value of $f(\thetahat)$ found in 5 runs of our algorithm. (When fitting to empirical final size distributions we find that the variability of $f(\thetahat)$ between runs is very small and thus 5 runs is more than sufficient to be confident that we have found the best-fitting model; see Section~\ref{sec:ident} for further details.) For clarity, the results in the figure have been ordered according to the best fit of the true model.

From these plots it is immediately clear that the correct model (i.e.\ the one that generated the data) has the best fit on most occasions (24 out of 25 for the IDS-HH data and 22 out of 25 for the MT-HH data). It also seems clear that it is generally easier to rule out the MT-HH model when looking at data from the IDS-HH model than vice-versa (the gap between the two lines is generally much larger in plot (a) than in plot (b)). Especially intriguing is the observation that, when the data are from the MT-HH model there is a clear association between the K-L distances to the best IDS-HH and MT-HH model; however there is much less, if any, association when the data are from the IDS-HH model. This may be an artifact of the fact that, as can be seen from Figure~\ref{fig:InfHist}, the IDS-HH model is generally able to fit MT-HH data better than the MT-HH model can fit IDS-HH data.

\subsection{Pseudolikelihood motivation for use of K-L divergence}
\label{sec:KLmotivation}
In this subsection we motivate our choice of the Kullback-Leibler divergence for assessing the distance between the two models by relating it to a maximum pseudolikelihood estimation procedure. For ease of presentation our arguments are informal, rather than fully rigorous.  Suppose, as above, that we have data $\{ q_n (r_M , r_S ) \}$ from an epidemic in a community of $m$ households. If we make the approximation that the outcomes in different households are mutually independent then the likelihood of these data under one or the other of our models is given by
\begin{equation} 
L(\theta)=\prod_{n=1}^{n_{\max}} \prod_{r_M,r_S} [ p_n (r_M,r_S | \theta)]^{m_n q_n(r_M,r_S)}, \label{4.1}
\end{equation}
where $m_n$ is the number of households of size $n$ in the community, and $(\theta,\, p_n(r_M,r_S | \theta))$ is either $(\theta^{(MT)},\, p^{(MT)}_n(r_M,r_S | \theta^{(MT)}))$ or $(\theta^{(IDS)},\, p^{(IDS)}_n(r_M, \newline r_S | \theta^{(IDS)}))$. In reality, \eqref{4.1} is a pseudolikelihood since the outcomes in distinct households are dependent, as they are part of the same community-wide epidemic, though the dependence is small (the covariance of the final outcomes in \emph{distinct} households is of order $1/m$ for large $m$; cf.\ Ball and Lyne (2010)).  The maximum pseudolikelihood estimator of $\theta$, denoted by $\thetahat$, is obtained by maximising $L(\theta)$, or equivalently by maximising $l(\theta)=\log L(\theta)$, which is given by
\[
l(\theta)=m \sum_{n=1}^{n_{\max}} \rho_n \sum_{r_M,r_S} q_n (r_M,r_S) \log p_n (r_M,r_S | \theta).
\]
Note that maximising $l(\theta)$ is equivalent to minimising the Kullback--Leibler divergence 
%JMB version
%$D_{KL} (q || p(\theta))$,
%Normal version
\newline $D_{KL} (q || p(\theta))$,
 defined by~\eqref{eq:KLdef}. Moreover, the pseudolikelihood ratio goodness-of-fit test statistic, $\Lambda_m$ say, for assessing the adequacy of the model for these data is given by
\begin{equation}
-2 \log \Lambda_m = 2m D_{KL} (q || p ( \thetahat)); \label{4.2}
\end{equation}
cf., for example, Bishop \etal\ (1975, Equation~10.2-6), who consider testing the goodness-of-fit of a specified multinomial model.

Now, for example, suppose that these data $\{ q_n (r_M,r_S)\}$ were actually generated by the IDS-HH model with parameter $\theta^{(IDS)}$, but that we fit the MT-HH model.  Then (cf.\ Section~\ref{sec:IDSLargePop}) $q_n (r_M,r_S) \overset{p}{\longrightarrow} p_n^{(IDS)} (r_M,r_S | \theta^{(IDS)})$ as $m \to \infty$, whence $\hat{\theta}^{(MT)} \overset{p}{\longrightarrow} \theta_\ast^{(MT)}$ as $m \to \infty$, where $\theta_\ast^{(MT)}$ minimises $D_{KL} ( p^{(IDS)} ( \theta^{(IDS)} ) || p^{(MT)} ( \theta^{(MT)}))$ with respect to $\theta^{(MT)}$.
%, where $D_{KL} (p^{(IDS)} ( \theta^{(IDS)}) || p^{(MT)} (\theta^{(MT)}))$ is obtained by setting $f_n (r_M,r_S)=p_n^{(IDS)} (r_M,r_S | \theta^{(IDS)})$ in $D_{KL} (f || p^{(MT)} (\theta^{(MT)}))$. 
(Here, $\overset{p}{\longrightarrow}$ denotes convergence in probability.)  Hence, using (\ref{4.2}),
\[
- \frac{1}{m} 2 \log \Lambda_m \overset{p}{\longrightarrow} 2D_{KL} (p^{(IDS)} ( \theta^{(IDS)}) || p^{(MT)} (\theta_{\ast}^{(MT)}))~~ \text{ as }~~ m \to \infty .
\]

If, instead, we fit the IDS-HH model with parameter $\theta^{(IDS)}$, then $\hat{\theta}^{(IDS)} \overset{p}{\longrightarrow} \theta^{(IDS)}$ as $m \to \infty$ and, since $D_{KL} (p^{(IDS)} (\theta^{(IDS)}) || p^{(IDS)} (\theta^{(IDS)}))=0$, $-2m^{-1} \newline \log \Lambda_m \overset{p}{\longrightarrow} 0$ as $m \to \infty$.
In these circumstances, $-2 \log \Lambda_m$ asymptotically equals the usual chi-square goodness of fit test statistic
\[
X^2 = \sum_{n=1}^{n_{\max}} \sum_{r_M,r_S} \frac{(m_n q_n(r_M,r_S)-m_n p_n^{(IDS)} (r_M,r_S | \hat{\theta}^{(IDS)}))^2} {m_n p_n^{(IDS)} (r_M,r_S | \hat{\theta}^{(IDS)})}.
\]
However, dependencies between the households imply that $X^2$ may not have the usual asymptotic $\chi^2$ distribution; instead the asymptotic distribution of $X^2$ is a linear combination of $d$ independent $\chi_1^2$ random variables, where $d$ is the degrees of freedom of the usual chi-square test (cf.\ Ball and Lyne (2010)).  Nevertheless, (\ref{4.2}) gives a guide for interpreting both our infinite and finite population model discrimination results.  Moreover, if these data $\{ q_n (r_M,r_S)\}$ come from a small fraction, $\epsilon_0$ say, of the households among which the epidemic is spreading, as is often the case in practice, then, if the model is correct, the asymptotic distribution of $X^2$ is very close to the usual $\chi_{d}^2$ distribution, the approximation being exact in the limit as $\epsilon_0\downarrow0$ (cf.\ Ball and Lyne (2010)).

Clearly there are precisely analagous results which hold if the data instead come from the MT-HH model.

\subsection{Identifiability and model fitting}
\label{sec:ident}

In this subsection we investigate how our model fitting methodology works in practice. We find that there are two key issues that influence the overall behaviour of our algorithm. The first is the striking difference in the distribution of final K-L distances $f(\thetahat)$ in the situations where the target final size distribution $q$ can or cannot be captured \emph{exactly} (to numerical accuracy) by the asymptotic version of the model we try to fit. The target final size distribution cannot be captured exactly when either (i) we try to fit the wrong model or (ii) the data is from a finite population. (This has the important consequence that when fitting a model to empirical final size distributions we need only run the algorithm a few, say 5--10, times to be confident that we have found the best possible fit.) We therefore restrict our attention here to an exploration of seeking to fit the models to `data' which are the asymptotic ($m\to\infty$) distributions corresponding to the parameter values given earlier, with household size distribution $\rho^{(3)}=(1,1,1)/3$. We refer to these data with this household size distribution and the parameters given previously as $q^{(MT)}$ and $q^{(IDS)}$. In the course of this we clearly see the second issue that arises, namely that there are identifiability issues in the IDS-HH model. In the IDS-HH model as we parameterise it, it seems that some parameters are identifiable while some are more difficult to identify, though we can find functions of these parameters which do appear to be identifiable.

\subsubsection{Fitting the correct model}
\begin{figure}
\begin{center}
\psfrag{log10fthetahat}[][]{\tiny $\log_{10} f(\thetahat)$}
\includegraphics[width=\hfigwidth]{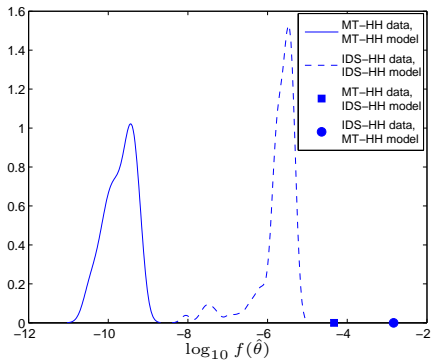}
\caption{Profiles of final K-L distances obtained when fitting both models to both $q^{(MT)}$ and $q^{(IDS)}$.}
\label{fig:exactFitHist}
\end{center}
\end{figure}

Firstly we look at fitting each model to data generated from that same model, so we should be able to recover the input parameters used to generate the data and find that $f(\thetahat)$ is very close to 0. Figure~\ref{fig:exactFitHist} shows density estimates (essentially smoothed histograms, which we use for ease of display) of $f(\thetahat)$ for the best 90 of 100 runs of our algorithm when fitting each model to data generated by that model. We use only the best 90\% of runs so as to exclude the poor fits sometimes obtained for the reasons explained in the second paragraph of Section~\ref{sec:modelFitting}. (For comparison, Figure~\ref{fig:exactFitHist} also shows the smallest $f(\thetahat)$ values found when fitting each model to the data generated by the \emph{other} model; these are displayed as points rather than densities since, as shown in Section~\ref{sec:FitWrongModel}, in these circumstances the variability of $f(\thetahat)$ is very small.)
This figure shows that our algorithm consistently finds model parameters $\thetahat$ which quite accurately reproduce the target final size distributions, indicated by the very small values of $f(\thetahat)$. We see shortly why the MT-HH model can be fitted to its own final size distribution rather better then the IDS-HH model. If we examine the parameter estimates $\thetahat$ that yield these final K-L distances we see (Table~\ref{tab:MTfitToMTInfData}) that the MT-HH model recovers the parameters used to generate $q^{(MT)}$ with a high degree of accuracy and very little variability, whereas when we fit the IDS-HH model to $q^{(IDS)}$ we find (Table~\ref{tab:IDSfitToIDSInfData}) that several parameters are estimated quite poorly.

\begin{table}[h]
\begin{center}
\caption{Summary of parameter estimates when fitting the MT-HH model to $q^{(MT)}$ (best 90 of 100 runs).}
\begin{equation*}
\begin{array}{r|ccccccc}
{\rm Parameter} &
\pi_M & \pi_S & \lambda^{(L)}_{MM} & \lambda^{(L)}_{MS} & \lambda^{(L)}_{SM} & \lambda^{(L)}_{SS} & \beta_M \\ \hline
{\rm True\ value} & 
0.7263 & 0.5224 & 0.2000 & 0.4000 & 0.4000 & 0.8000 & 0.4000 \\
{\rm Mean} & 
0.7263 & 0.5224 & 0.2000 & 0.4000 & 0.4000 & 0.8000 & 0.4000 \\
{\rm Std.\ dev.} & 
0.00003 & 0.00004 & 0.00004 & 0.00013 & 0.00008 & 0.00022 & 0.00004
\end{array}
\end{equation*}
\label{tab:MTfitToMTInfData}
\end{center}
\end{table}
\begin{table}[h]
\begin{center}
\caption{Summary of parameter estimates when fitting the IDS-HH model to $q^{(IDS)}$ (best 90 of 100 runs).}
\begin{equation*}
\begin{array}{r|ccccccccc}
{\rm Parameter} &
\lambda^{(G)}_M & \lambda^{(G)}_S & \lambda^{(L)}_M & \lambda^{(L)}_S & p^{(G)}_{MM} & p^{(G)}_{SM} & p^{(L)}_{MM} & p^{(L)}_{SM} & \gamma_S \\ \hline
{\rm True\ value} & 
1.0000 & 2.0000 & 0.5000 & 1.0000 & 0.8000 & 0.2000 & 0.5000 & 0.1000 & 1.5000 \\
{\rm Mean} & 
1.7805 & 3.9959 & 0.5028 & 4.2510 & 0.3924 & 0.4954 & 0.4935 & 0.0932 & 8.5436 \\
{\rm Std.\ dev.} & 
0.6781 & 3.4214 & 0.0017 & 2.4507 & 0.1662 & 0.2507 & 0.0040 & 0.0046 & 4.9344 
\end{array}
\end{equation*}
\label{tab:IDSfitToIDSInfData}
\end{center}
\end{table}

\subsubsection{Identifiability in the IDS-HH model}
The poorer recovery of input parameters in the IDS-HH model can to a large extent be explained by issues of identifiability. We mention in Section~\ref{sec-MT-HH} that it is known that in the MT-HH model the global rate parameters are not uniquely identifiable from final size data (Ball \etal\ 2004) and for this reason we estimate the probabilities $\bpi$ rather than the global rates $(\lambda^{(G)}_{MM}, \lambda^{(G)}_{MS}, \lambda^{(G)}_{SM}, \lambda^{(G)}_{SS})$. However, the IDS-HH model we propose is new so such identifiability issues have not been explored. Moreover, identifiability is difficult to study rigorously for this model as there is no analytical expression for the household final size distributions. Careful examination of the parameter estimates when fitting the IDS-HH model to $q^{(IDS)}$ suggests that some identifiability issues are present here. In particular, in our parameterisation of the IDS-HH model there are three combinations of parameters that seem identifiable whilst some of the individual parameters are very difficult to identify separately.

%\emph{Note also that $p^{(G)}$'s are estimated poorly. Unidentifiablity similar to the global parameters in the MT-HH model? Variation in $p^{(G)}$'s explainable by variation in $\lambda^{(G)}$'s which is in turn explainable by variation in $\gamma_S$?? Investigate and discuss briefly; for a start, look at pairwise plots of $\lambda^{(G)}$'s and $p^{(G)}$'s. Even fixing $\gamma_S$, which helps enormously with the estimation of $\lambda^{(L)}_S$ and the linear combination of $\lambda^{(G)}$'s, has little impact on the $p^{(G)}$'s.}

The first of these combinations is $\lambda^{(L)}_S$ and $\gamma_S$; our algorithm estimates the ratio $\lambda^{(L)}_S/\gamma_S$ extremely well (see Table~\ref{tab:IDSfitToIDSInfData2}), but has difficulty identifying the precise values of these parameters. The second set of troublesome parameters consists of the global contact rates $\lambda^{(G)}_M$ and $\lambda^{(G)}_S$ and the removal rate $\gamma_S$. If the removal rates $\gamma_M$ and $\gamma_S$ are known then the relationship $\pi_G = \exp(-(z_M \lambda^{(G)}_M/\gamma_M + z_S \lambda^{(G)}_S/\gamma_S))$, where $\pi_G = q_1(0,0)$ ($=\sqrt[n]{q_n(0,0)}$ for any $n\leq \nmax$) is the probability that a given individual avoids global infection, specifies a linear equation that $\lambda^{(G)}_M$ and $\lambda^{(G)}_S$ must satisfy. If we assume that the removal rates are both known then our algorithm identifies the correct linear combination of global contact rates very easily but finds it very difficult, though possible, to find the most likely values of these parameters individually. However, we assume that (one of) the removal rates is unknown and, as just discussed, not estimated very well; thus the global rates are generally not estimated very reliably. Nevertheless, when the initial guess for $\gamma_S$ is close to its optimum (correct) value, $\lambda^{(L)}_S$ and the above linear combination of $\lambda^{(G)}_M$ and $\lambda^{(G)}_S$ are also estimated easily and reasonably accurately, and a very good fit is obtained. In the latter situation it is also possible to recover the individual rates $\lambda^{(G)}_M$ and $\lambda^{(G)}_S$ with our algorithm but this is far more difficult. That $z_M \lambda^{(G)}_M/\gamma_M + z_S \lambda^{(G)}_S/\gamma_S$ is estimated well is demonstrated in Table~\ref{tab:IDSfitToIDSInfData2}, in which $z_M$ and $z_S$ are given by their observed values in the (infinite) data.

The other parameters with identifiability issues are the global infection probabilities $p^{(G)}_{MM}$ and $p^{(G)}_{SM}$. Some information about these parameters can be obtained by considering households of size 1 which become infected. Focus on such a household and suppose that there are in total $Y_M$ mild and $Y_S$ severe infectives in the population just prior to its infection. Then the probability that this infection is mild is
\begin{equation*}
\frac{Y_M \lambda^{(G)}_M p^{(G)}_{MM} + Y_S \lambda^{(G)}_S p^{(G)}_{SM}}{Y_M \lambda^{(G)}_M  + Y_S \lambda^{(G)}_S}.
\end{equation*}
Now, $Y_M$ and $Y_S$ are random and vary throughout the epidemic. A crude approximation is to replace the ratio $Y_M/Y_S$ by $\gamma_M^{-1}z_M/\gamma_S^{-1}z_S$, the latter taking into account the different infectious periods of mild and severe infectives. Thus the proportion of infected households of size 1 that are mildly infected is approximately
\begin{equation*}
\frac{z_M \lambda^{(G)}_M p^{(G)}_{MM}/\gamma_M + z_S \lambda^{(G)}_S p^{(G)}_{SM}/\gamma_S}{z_M \lambda^{(G)}_M/\gamma_M  + z_S \lambda^{(G)}_S/\gamma_S},
\end{equation*}
leading to the relationship
\begin{equation}
\frac{p_1^{(IDS)}(1,0 | \theta^{(IDS)})}{p_1^{(IDS)}(1,0 | \theta^{(IDS)}) + p_1^{(IDS)}(0,1 | \theta^{(IDS)})}
\approx
\frac{z_M \lambda^{(G)}_M p^{(G)}_{MM}/\gamma_M + z_S \lambda^{(G)}_S p^{(G)}_{SM}/\gamma_S}{z_M \lambda^{(G)}_M/\gamma_M  + z_S \lambda^{(G)}_S/\gamma_S}.
\label{eq:pGs}
\end{equation}
We have seen above that the denominator in the right hand side of~\eqref{eq:pGs} can be estimated well, hence it is reasonable to expect that the numerator might be too. That this is indeed the case is borne out in Table~\ref{tab:IDSfitToIDSInfData2}.

\begin{table}[h]
\begin{center}
\caption{Functions of estimated IDS-HH model parameters when fitting IDS-HH model to $q^{(IDS)}$ (best 90 of 100 runs).}
\begin{equation*}
\begin{array}{c|c|c|c}
{\rm Function} &
z_M\lambda^{(G)}_M/\gamma_M + z_S\lambda^{(G)}_S/\gamma_S & \lambda^{(L)}_S/\gamma_S & z_M \lambda^{(G)}_M p^{(G)}_{MM} / \gamma_M + z_S\lambda^{(G)}_S p^{(G)}_{SM} / \gamma_S\\ \hline
{\rm True\ value} & 
0.50669 & 0.50000 & 0.21340 \\
{\rm Mean} & 
0.50672 & 0.49807 & 0.21069 \\
{\rm Std.\ dev.} & 
0.00003 & 0.00113 & 0.00170
\end{array}
\end{equation*}
\label{tab:IDSfitToIDSInfData2}
\end{center}
\end{table}

\subsubsection{Fitting the incorrect model}
\label{sec:FitWrongModel}
We now turn our attention to the situation where we try to fit one of the models to final size data arising from the other model. Fitting the MT-HH model to IDS-HH data gives parameter estimates summarised in Table~\ref{tab:MTfitToIDSInfData}. While there is more variation in the MT-HH parameter estimates than when we fit to data from the MT-HH model, the variation is still relatively small. Furthermore, the variation in the final K-L distances $f(\thetahat)$ is very small (mean $1.46 \times 10^{-3}$, st.\ dev.\ $1.5 \times 10^{-10}$), giving confidence that (i) we have found the region of parameter space where the MT-HH model can best reproduce the data from the IDS-HH model and (ii) that the best fitting MT-HH model does not reproduce the IDS-HH final size distribution very well. For comparison the minimum of these K-L distances is also shown in Figure~\ref{fig:exactFitHist}, as a point rather than a density because the variation is so small.
\begin{table}[h]
\begin{center}
\caption{Summary of parameter estimates when fitting the MT-HH model to $q^{(IDS)}$ (best 90 of 100 runs).}
\begin{equation*}
\begin{array}{r|ccccccc}
{\rm Parameter} &
\pi_M & \pi_S & \lambda^{(L)}_{MM} & \lambda^{(L)}_{MS} & \lambda^{(L)}_{SM} & \lambda^{(L)}_{SS} & \beta_M \\ \hline
{\rm Mean} & 
0.5210 & 0.6450 & 1.3712 & 0.2561 & 0.0509 & 0.8990 & 0.3373 \\
{\rm Std.\ dev.} & 
0.00003 & 0.00001 & 0.00035 & 0.00003 & 0.00003 & 0.00007 & 0.00002
\end{array}
\end{equation*}
\label{tab:MTfitToIDSInfData}
\end{center}
\end{table}

Lastly we consider fitting the IDS-HH model to the MT-HH data; see Table~\ref{tab:IDSfitToMTInfData} and Figure~\ref{fig:exactFitHist}. Here we see variations in $\thetahat$ roughly the same as those seen when fitting the IDS-HH model to data it can reproduce exactly. Again we find that the estimates of $\lambda^{(L)}_M$ and the $p^{(L)}$'s show little variation and we also see the same identifiability issues present. Although the estimates of $\gamma_S$, $\lambda^{(L)}_S$, $\lambda^{(G)}_S$ and the $p^{(G)}$'s individually vary wildly we find that $\lambda^{(L)}_S/\gamma_S$, $z_M \lambda^{(G)}/\gamma_M + z_S \lambda^{(G)}_S/\gamma_S$ and $z_M \lambda^{(G)}_M p^{(G)}_{MM}/\gamma_M + z_S \lambda^{(G)}_S p^{(G)}_{SM}/\gamma_S$ show very little variation (see Table~\ref{tab:IDSfitToMTInfData2}). Similarly to when we fit the MT-HH model to the IDS-HH data, we find very little variability in the final K-L distances that we find, for the 90 smallest values the mean and st.\ dev.\ are $4.69\times 10^{-5}$ and $1.0\times10^{-7}$, respectively. Though these K-L distances certainly seem bounded away from zero, suggesting that the IDS-HH model cannot reproduce the MT-HH final size distribution, they are appreciably smaller than when fitting the MT-HH model to IDS-HH data.%, though this is not surprising as the IDS-HH model has more parameters that we fit.

\begin{table}[h]
\begin{center}
\caption{Summary of parameter estimates when fitting the IDS-HH model to $q^{(MT)}$ (best 90 of 100 runs).}
\begin{equation*}
\begin{array}{r|ccccccccc}
{\rm Parameter} &
\lambda^{(G)}_M & \lambda^{(G)}_S & \lambda^{(L)}_M & \lambda^{(L)}_S & p^{(G)}_{MM} & p^{(G)}_{SM} & p^{(L)}_{MM} & p^{(L)}_{SM} & \gamma_S \\ \hline
{\rm Mean} & 
2.3496 & 4.0997 & 0.1982 & 4.2774 & 0.2351 & 0.5028 & 0.2757 & 0.3470 & 9.7983 \\
{\rm Std.\ dev.} & 
0.6919 & 2.2107 & 0.0001 & 2.8475 & 0.1548 & 0.2579 & 0.0038 & 0.0037 & 6.5202
\end{array}
\end{equation*}
\label{tab:IDSfitToMTInfData}
\end{center}
\end{table}

\begin{table}[h]
\begin{center}
\caption{Functions of estimated IDS-HH model parameters when fitting IDS-HH model to $q^{(MT)}$ (best 90 of 100 runs).}
\begin{equation*}
\begin{array}{c|c|c|c}
{\rm Function} &
z_M\lambda^{(G)}_M/\gamma_M + z_S\lambda^{(G)}_S/\gamma_S & \lambda^{(L)}_S/\gamma_S & z_M \lambda^{(G)}_M p^{(G)}_{MM} / \gamma_M + z_S\lambda^{(G)}_S p^{(G)}_{SM} / \gamma_S \\ \hline
{\rm Mean} & 
0.50504 & 0.57068 & 0.13909 \\
{\rm Std.\ dev.} & 
0.00003 & 0.00023 & 0.00143
\end{array}
\end{equation*}
\label{tab:IDSfitToMTInfData2}
\end{center}
\end{table}

%\emph{Perhaps a further comment that the behaviour (variability) of both $\thetahat$ and $f(\thetahat)$ when fitting the models to finite data is essentially the same as when fitting to the wrong model since the empirical final size distributions (even when fitting the correct model) are not \emph{exactly} reproducible with the model we fit.}

\section{Discussion}
\label{disc}

In this paper we define two candidate models that might explain how an infectious disease having varying disease response could spread in a community of households. Large population properties of the two models are presented. These results are used to show by means of numerical illustrations, that it is generally possible to discriminate between the two models. More precisely, given final outcome data from a sufficiently large community of households it is, except in some degenerate cases, possible to determine which of the two explanations to varying disease response that best explain the data.

Both models could of course be extended towards higher realism. For example, besides household structure, all individuals are assumed similar whereas it would be more realistic to distinguish between adults and children having different mixing rates. Another extension would be to allow for more than two different disease responses. It is of course also possible to come up with other models giving rise to mild and severe infectives. However, we believe that the two models studied capture the perhaps two most likely reasons: either the infection status of an individual is predetermined or else it depends on whom the person was infected by. In the first situation it could be natural to extend the model to allow this predetermined status to be dependent between individuals of the same household, for example due to previous exposure to the disease. In the present model it is assumed that the predetermined infection status is independent also within households. Another important extension would of course be to apply the method to real data with the hope to find out more about the underlying reason for having varying disease response.

\section*{Acknowledgements}
%\begin{acknowledgements}
This research was supported by the Swedish Research Council and by the UK
Engineering and Physical Sciences Research Council (under grants
EP/F03234X/1 and EP/E038670/1).
%\end{acknowledgements}

\section*{References}
%\begin{thebibliography}{99}
\indent

%\bibitem{Addy91}
Addy C.~L.,\ Longini I.~M.\ and Haber M.\ (1991).
A generalized stochastic model for the analysis of infectious disease final size data.
{\it Biometrics} 47: 961--74.

%\bibitem{BB06}
Ball F.~G.\ and Becker, N.~G.\ (2006).
Control of transmission with two types of infection.
{\it  Math. Biosci.} 200:170--187.

%\bibitem{BB07}
Ball F.~G.\ and Britton, T.\ (2007).
An epidemic model with infector-dependent severity.
{\it  Adv. Appl. Prob.} 39:949--972.

%\bibitem{BB09}
Ball F.~G.\ and Britton, T.\ (2009).
An epidemic model with infector and exposure dependent severity.
{\it  Math. Biosci.} 218:105--120.

%\bibitem{BBL04}
Ball F.~G., Britton, T.\ and Lyne, O.~D.\ (2004).
Stochastic multitype epidemics in a community of households: Estimation of threshold parameter $R_*$ and secure vaccination coverage.
{\em Biometrika} 91:345--362.

%\bibitem{BL01}
Ball, F.~G.\ and Lyne, O.~D.\ (2001).
Stochastic multitype SIR epidemics among a population partitioned into households.
{\em Adv. Appl. Prob.} 33:99--123.

%\bibitem{BL10}
Ball, F.~G.\ and Lyne, O.~D.\ (2010).
Statistical inference for epidemics among a population of households.
Under revision.

%\bibitem{BH96}
Becker, N.~G.\ and Hall, R.\ (1996).
Immunization levels for preventing epidemics in a community of households made up of individuals of various types.
{\em Math. Biosci.} 132:205--216.

%\bibitem{BFH75}
Bishop, Y.~M.~M., Feinberg, S.~E.\ and Holland, P.~W.\ (1975).
{\em Discrete multivariate statistics}.  MIT Press, Cambridge.

%\bibitem{BB00}
Britton, T.\ and Becker, N.~G.\ (2000).
Estimating the immunity coverage required to prevent epidemics in a community of households.
{\em Biostatistics} 1:389--402.

%\bibitem{CF08}
Carrat, F., Vergu, E., Ferguson, N.~M., Lemaitre, M., Cauchemez, S., Leach, S.\ and Valleron, A.-J.\ (2008).
Time lines of infection and disease in human influenza: a review of volunteer challenge studies.
{\em Am.\ J.\ Epidemiology} 167:775--785.

%\bibitem{EK86}
Ethier, S.~N.\ and Kurtz, T.~G.\ (1986).
{\em Markov Processes: Characterization and Convergence}. Wiley, New York.

%\bibitem{Fetal05}
Ferguson, N.~M., Cummings, D.~A.~T., Cauchemez, S., Fraser, C., Riley, S., Meeyai, A., Iamsirithaworn, S., Burke, D.~S.\ (2005).
Strategies for containing an emerging influenza pandemic in Southeast Asia.
{\em Nature} 437:209--214.

%\bibitem{LBDLM01}
Leroy, E.~M., Baize, S., Debre, P., Lansoud-Soukate, J.\ and Mavoungou, E.\ (2001).
Early immune responses accompanying human asymptomatic Ebola infections.
{\em Clinical and Experimental Immunology} 124:453--460.

%\bibitem{MI98}
Mangada, M.~N.~M.\ and Igarashi, A.\ (1998).
Molecular and in vitro analysis of eight dengue type 2 viruses isolated from patients exhibiting different disease severities.
{\em Virology} 244:458--466.

Mehta, P.~N.\ and Chatterjee, A.\ (2010).
Varicella. {\em eMedicine\/}. Accessed 29th April 2010.
http://emedicine.medscape.com/article/969773-overview

%\bibitem{MA97}
Morley, D.~C.\ and Aaby, P.\ (1997).
Managing measles: size of infecting dose may be important.
{\em BMJ\/} 314:1692.

%\bibitem{PFF08}
Pellis, L.,\ Ferguson, N.~M.\ and Fraser, C.\ (2008).
 The relationship between real-time and discrete-generation models of epidemic spread.
{\it Math. Biosci.} 216:63--70.

%\bibitem{SH98}
Staalsoe, T.\ and Hviid, L.\ (1998).
The Role of Variant-specific Immunity in Asymptomatic Malaria Infections: Maintaining a Fine Balance.
{\em Parasitology Today} 14:177--178.

%\end{thebibliography}

\end{document}